\input ./style/arxiv-general.cfg
\documentclass[labelcite,seceqn,MSNbibl,number,citesort,dvips]{arxbj}
\makeatletter
   \@ifpackageloaded{graphicx}{}{\usepackage{graphicx}}
\makeatother
\usepackage{upgreek}
\usepackage{arydshln}

\aid{0}
\volume{21}
\issue{4}
\pubyear{2015}
\firstpage{2093}
\lastpage{2119}
\doi{10.3150/14-BEJ636} 
\docsubty{FLA}

\makeatletter
\newcommand{\rrvert}{\vert}
\newcommand{\rrVert}{\Vert}
\newcommand{\llvert}{\vert}
\newcommand{\llVert}{\Vert}

\def\mod{\operatorname{mod}}
\def\bmod{\quad\operatorname{mod}}

\newcommand{\prolim}{p\mbox{-}\lim}

\newcommand{\field}[1]{\mathbb{#1}}
\newcommand{\R}{\field{R}}
\newcommand{\N}{\field{N}}
\newcommand{\Z}{\field{Z}}

\newcommand{\tr}{\operatorname{tr}}
\newcommand{\var}{\operatorname{Var}}
\newcommand{\cov}{\operatorname{Cov}}
\newcommand{\cum}{\operatorname{cum}}
\newcommand{\Int}{\operatorname{Int}}
\newcommand{\E}{\mathrm{E}}

\renewcommand{\bm}[1]{\boldsymbol{#1}}
\newcommand{\bmm}[1]{\mathbf{#1}}

\newcommand{\xrightarrow}[1]{\stackrel{#1}{\rightarrow}}
\newcommand{\plim}{\xrightarrow{\mathcal{P}}}
\newcommand{\dlim}{\xrightarrow{\mathcal{L}}}

\newproclaim{ass}{Assumption}[section]
\newtheorem{lem}{Lemma}[section]

\newtheorem{theorem}{Theorem}[section]

\newremark{rem}{Remark}[section]

\newcommand{\eqref}[1]{(\ref{#1})}

\renewcommand{\epsilon}{\varepsilon}
\makeatother

\begin{document}
\begin{frontmatter}

\title{An empirical likelihood approach for symmetric $\alpha$-stable processes}
\runtitle{Empirical likelihood approach for $\alpha$-stable processes}

\begin{aug}
\author[A]{\inits{F.}\fnms{Fumiya}~\snm{Akashi}\corref{}\thanksref{e1}\ead[label=e1,mark]{fakashi01@fuji.waseda.jp}},
\author[A]{\inits{Y.}\fnms{Yan}~\snm{Liu}\thanksref{e2}\ead[label=e2,mark]{y.liu2@kurenai.waseda.jp}} \and
\author[A]{\inits{M.}\fnms{Masanobu}~\snm{Taniguchi}\thanksref{e3}\ead[label=e3,mark]{taniguchi@waseda.jp}}
\address[A]{Department of Applied Mathematics, School of Fundamental
Science and Engineering,
Waseda University 3-4-1, Okubo, Shinjuku-ku, Tokyo, 169-8555, Japan.\\
\printead{e1,e2,e3}}
\end{aug}

\received{\smonth{4} \syear{2013}}
\revised{\smonth{4} \syear{2014}}

%
\begin{abstract}
Empirical likelihood approach is one of non-parametric statistical methods,
which is applied to the hypothesis testing or construction of
confidence regions for pivotal unknown quantities.
This method has been applied to the case of independent identically
distributed random variables
and second order stationary processes.
In recent years, we observe heavy-tailed data in many fields.
To model such data suitably,
we consider symmetric scalar and multivariate $\alpha$-stable linear
processes generated
by infinite variance innovation sequence.
We use a Whittle likelihood type estimating function in the empirical
likelihood ratio function and
derive the asymptotic distribution of the empirical likelihood ratio
statistic for $\alpha$-stable linear processes.
With the empirical likelihood statistic approach,
the theory of estimation and testing for second order stationary processes
is nicely extended to heavy-tailed data analyses, not straightforward,
and applicable to a lot of financial statistical analyses.
\end{abstract}

%
\begin{keyword}
\kwd{confidence region}
\kwd{empirical likelihood ratio}
\kwd{heavy tail}
\kwd{normalized power transfer function}
\kwd{self-normalized periodogram}
\kwd{symmetric $\alpha$-stable process}
\kwd{Whittle likelihood}
\end{keyword}
\end{frontmatter}

\section{Introduction}\label{sec1}
Non-parametric methods have been developed for the
statistical analysis of univariate and multivariate observations in the
area of time series analysis
to carry out the problem of inference and hypothesis testing.
Rank-based methods and
empirical likelihood methods
have been introduced in succession in these two decades.

Owen \cite{O:1988} introduced the empirical likelihood approach to independent
and identically distributed ({i.i.d.}) data
and he showed that the empirical likelihood ratio statistic is
asymptotically $\chi^2$-distributed.
For dependent data, Monti \cite{M:1997}, Ogata and Taniguchi \cite{OT:2010} derived
the limit distribution of the empirical likelihood ratio statistic
based on the derivative of the Whittle likelihood with respect to parameters.
From these papers, we can construct confidence sets for the
coefficients in a predictor and autocorrelation coefficients in
multivariate stationary processes, etc.

In the last few decades,
heavy-tailed data have been observed
in a variety of fields involving electrical engineering, hydrology,
finance and physical systems (Nolan \cite{N:2012} and Samorodnitsky and Taqqu
\cite{samoradnitsky1994stable}).
In particular, Fama \cite{Fama1965} and Mandelbrot \cite{M1963} gave economic and
financial examples that
show such data are poorly grasped by Gaussian model.
When we fit a GARCH-model to some financial data
and estimate the stable index of the residuals by Hill's estimator
$\hat{\alpha}$,
we often observe that the tail of the distribution is heavier than that
of Gaussian model.

\begin{figure}

\includegraphics{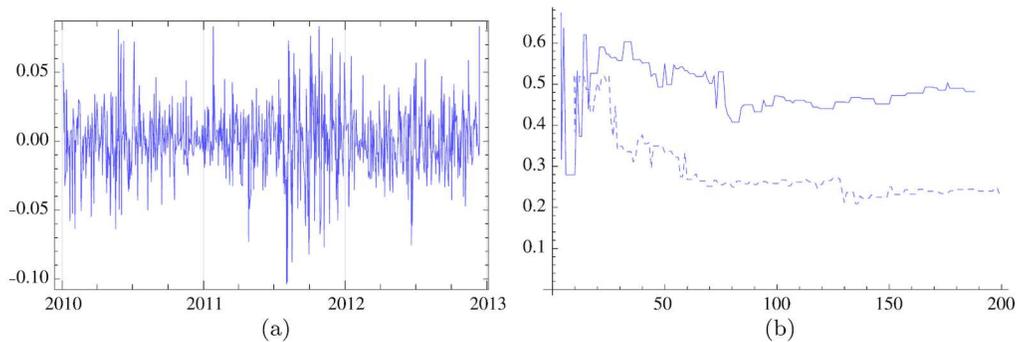}

\caption{Log return of Hewlett Packard company and the
Hill-plot. (a) log return of Hewlett Packard's stock price (from 1,
January, 2010 to 14, December, 2012).  (b) Hill-plot for residuals (dashed line is for i.i.d. normal
random variables).}\label{fig:01}
\end{figure}

Figure~\ref{fig:01} shows daily stock returns of Hewlett Packard company and the
Hill-plot for the residuals
(we used AIC to select the order of GARCH).
These graphs imply that it is more suitable to suppose these data are
generated from a process with stable innovations rather than to assume
these data have finite variances (for discussion of Hill-plot, see
Drees, de Haan and Resnick \cite{DHR:2000}, Hall \cite{label2568}, Hsing \cite{H:1991},
Resnick and St\v{a}rik\v{a} \cite{RS:1996} and \cite{RS:1998}).

To model such heavy-tailed data suitably, we introduce the following
linear process generated by stable innovations,
%
\begin{equation}\label{eq:1.1}
X(t) = \sum^{\infty}_{j=0}
\psi_j Z(t-j),\qquad t\in\Z,
\end{equation}
where $\psi_0 = 1$ and
$\{Z(t);\ t\in\Z\}$ ($\Z$ is the set of all integers) is a sequence of
i.i.d. symmetric $\alpha$-stable random variables
(for short s$\alpha$s).
In the case of $\alpha=2$, this process is Gaussian.
When $\alpha$ is less than 2, the usual spectral density function of
(\ref{eq:1.1}) cannot be defined.

Davis and Resnick \cite{DR1985a,DR1985b} and \cite{DR1986}
investigated the sample autocorrelation function (ACF) at lag $h$, and
derived the consistency of ACF.
Resnick and St{\u{a}}ric{\u{a}} \cite{RS:1998} gave a consistent
estimator of the tail index $\alpha$.
In view of the frequency domain approach,
Kl{\"u}ppelberg and Mikosch \cite{KM:1993,KM:1994} and \cite{KM:1996}
proposed a self-normalized periodogram because the expectation of the
usual periodogram does not exist, and introduced
some methods for parameter estimation and hypothesis testing.
Then, they showed that for any frequencies, self-normalized periodogram
converges to a random variable with finite second moment, and
proved the convergence of the functional of the self-normalized periodogram.

In this paper, we apply non-parametric method to the discrete linear
process \eqref{eq:1.1}.
It is natural to express the process non-parametrically
partly because finite parametric models often cannot describe real
data sufficiently, and partly because
there is no general solution of probability density function for
stable distribution.
Recently economists and quantitative analysts have introduced stable
stochastic models to asset returns in econometrics and finance.
In such situations, what we are interested in is to test statistical hypothesis
on an important pivotal quantity ``$\bm{\theta} = \bm{\theta}_0$'',
such as the correlation between the different realizations.
To achieve this goal,
Monti \cite{M:1997} and
Ogata and Taniguchi \cite{OT:2010}
employed the empirical likelihood to construct confidence sets for
linear processes
when innovations have finite variance.
A plausible way to define the important index $\bm{\theta}_0$ is
Whittle's approach, that is, $\bm{\theta}_0$ minimizes the disparity
%
\begin{equation}
\label{eq:disparity}
D({f}_{\bm{\theta}}, \tilde{g}) = \int_{-\uppi}^{\uppi}
\frac{\tilde{g}(\omega)}{{f}(\omega; \bm
{\theta})}\,\mathrm{d}\omega,
\end{equation}
where $\tilde{g}(\omega)$ is called a normalized power transfer
function of \eqref{eq:1.1},
and $f(\omega;\bm{\theta})$ is an appropriate score function.

This setting is useful for many situations.
For example, let us consider the $h$-step linear prediction of a scalar
stationary process $\{X(t); t\in\Z\}$.
We predict $X(t)$ by a linear combination of $\{X(s); s\leq t-h\}$,
\[
\hat X(t) = \sum^{\infty}_{j=h}
\phi_j(\bm{\theta}) X(t-j).
\]
The spectral representations of $X(t)$ and $\hat X(t)$ are
\[
X(t) = \int^{\uppi}_{-\uppi}\exp(-\mathrm{i}t\omega)\,\mathrm{d}
\zeta_X(\omega),\qquad \hat X(t) = \int^{\uppi}_{-\uppi}
\exp(-\mathrm{i}t\omega) \sum^{\infty
}_{j=h}
\phi_j(\bm {\theta})\exp(\mathrm{i}j\omega)\,\mathrm{d}\zeta_X(\omega),
\]
where $\{\zeta_X(\omega); -\uppi\leq\omega\leq\uppi\}$ is an orthogonal
increment process satisfying
\[
\E\,\mathrm{d}\zeta_X(\omega)\,\mathrm{d}\zeta_X(\mu) =
\cases{
g(\omega)\,\mathrm{d}\omega&\quad $(\omega= \mu)$,
\cr
0 &\quad $(\omega\neq\mu)$. } %
\]
Then, the prediction error is
%
\begin{equation}
\label{eq:2.5'} \E \bigl\llvert X(t)-\hat X(t) \bigr\rrvert ^2 = \int
^{\uppi}_{-\uppi} \Biggl\llvert 1-\sum
^{\infty
}_{j=h}\phi _j(\bm {\theta})\exp(\mathrm{i}j
\omega) \Biggr\rrvert ^2 g(\omega)\,\mathrm{d}\omega.
\end{equation}
Hence the best $h$-step predictor is given by $\sum^{\infty
}_{j=h}\phi
_j(\bm{\theta}_0)X(t-j)$, where
$\bm{\theta}_0$ minimizes (\ref{eq:2.5'}).
Comparing this with \eqref{eq:disparity}, if we set
\[
f(\omega;\bm{\theta}) = \Biggl\llvert 1-\sum^{\infty}_{j=h}
\phi_j(\bm {\theta })\exp(\mathrm{i}j\omega) \Biggr\rrvert ^{-2},
\]
this problem is exactly the same as that of seeking
$\bm{\theta}_0$ in their definition.
In addition to the linear prediction,
the empirical likelihood approach can also be applied to the case of
sample autocorrelation estimation,
which will be given in Section~\ref{sec2}.

The empirical likelihood ratio function for the problem of testing $H$:
$\bm{\theta}=\bm{\theta}_0$ is defined as
\[
R(\bm{\theta}) = \max_{w_1,\dots, w_n} \Biggl\{ \prod
_{t=1}^n n w_t; \sum
_{t=1}^n w_t \bmm{m}(
\lambda_t; \bm {\theta}) = \bm{0}, \sum_{t=1}^n
w_t = 1, 0 \leq w_1, w_2, \dots,
w_n \leq1 \Biggr\},
\]
and then the estimating function takes the form
\[
\bmm{m}(\lambda_t; \bm{\theta}) \equiv\frac{\partial}{\partial
\bm
{\theta}}
\frac{\tilde{I}_{n, X}(\lambda_t) }{{f}(\lambda_t; \bm{\theta})}, \qquad\lambda_t = \frac{2\uppi t}{n} \in(-\uppi, \uppi],
\]
where $\tilde{I}_{n,X}(\omega)$ is called self-normalized periodogram.
For our general stable linear process (\ref{eq:1.1}), we derive the
limit distribution of $R(\bm{\theta}_0)$
with its normalizing factor and construct the confidence interval
through a numerical method.

Here it should be noted that our extension to the stable case from the
finite variance case
requires new asymptotic methods,
and we report new aspects of the asymptotics for empirical likelihood approach,
which are different from the usual ones.
Furthermore, we extend the results to those of the mutivariate one with
independent innovations.
This is extremely important from a viewpoint of practical use.
In particular, we can analyze the relationship between two heavy-tailed
processes.
The way to derive the asymptotics of the multivariate case has also new aspects.
We find that the asymptotics for multivariate process need more
stronger conditions than
what we need in the univariate case.
The self-normalizing factor is also difficult to find in that case and we
use the norm of the stable series defined in Section~\ref{sec4} instead of the
square root matrix.

This paper is organized as follows:
In Section~\ref{sec2}, we shall introduce the fundamental setting and a brief
overview on
the empirical likelihood approach based on the Whittle likelihood.
With a different normalizing order for the empirical likelihood ratio function,
the main theoretical results, limit distribution of the empirical
likelihood ratio statistic for univariate and multivariate stable
linear processes, are formulated in Sections~\ref{sec3} and \ref{sec4}, respectively.
In Section~\ref{sec5}, the numerical results will be given under several settings.
We shall demonstrate some effectiveness of the empirical likelihood
ratio method.
The proofs of theorems in Sections~\ref{sec3} and \ref{sec4} are relegated to Section~\ref{sec6}.

As for notations and symbols used in this paper,
the set of all integers, non-negative integers ($=\{0,1,2,\ldots\}$) and
real numbers are denoted by $\Z$, $\N$ and $\R$, respectively.
For any sequence of random vectors $\{\bmm{A}(t); t\in\Z\}$,
$\mathbf{A}(t)
\plim\mathbf{A}$
and $\bmm{A}(t)\dlim\bmm{A}$, respectively, denote the convergence
to a random (or constant) vector $\bmm{A}$ in probability and law.
Especially, ``$\prolim_{t\to\infty}\bmm{A}(t) = \bmm{A}$''
implies ``$\bmm{A}(t)\plim\bmm{A}$ as $t\to\infty$''.
The transpose and conjugate transpose of matrix $\bmm{M}$ are denoted
by $\bmm{M}'$ and $\bmm{M}^*$,
and define $\Vert\bmm{M}\Vert_E:=\sqrt{\tr[\bmm{M}^*\bmm{M}]}$.

\section{Fundamental setting}\label{sec2}
In this section, we state the fundamental setting for the main results.
Throughout this paper, we use the following notations. For any sequence
$\{A(t); t\in\Z\}$ of random variables,
%
\begin{eqnarray}\label{eq:2.1}
\gamma_{n,A}^2 &=& n^{-2/\alpha}\sum
^{n}_{t=1}A(t)^2,
\nonumber
\\
I_{n,A}(\omega) &=& n^{-2/\alpha} \Biggl\llvert \sum
^{n}_{t=1}A(t) \exp (\mathrm{i}t\omega ) \Biggr\rrvert
^2,
\\
\tilde{A}_t &= &\frac{A(t)}{\sqrt{A(1)^2+\cdots+A(n)^2}}, \qquad t = 1,\ldots,n,\nonumber
\end{eqnarray}
and
\[
\tilde{I}_{n,A}(\omega) = \frac{I_{n,A}(\omega)}{\gamma_{n,A}^2} = \Biggl\llvert \sum
^{n}_{t=1}\tilde{A}_t \exp(\mathrm{i}t
\omega) \Biggr\rrvert ^2.
\]
We call $\tilde{I}_{n,A}(\omega)$ a \textit{self-normalized periodogram}
of $A(1),\ldots,A(n)$.
Mikosch \textit{et al.} \cite{MGKA:1995} studied estimation of the following stable and
causal ARMA process:
\begin{eqnarray*}
&&X(t) + \phi_1 X(t-1) + \cdots+ \phi_p X(t-p) = Z(t) +
\theta_1 Z(t-1) + \cdots+ \theta_q Z(t-q),
\\
&&\quad Z(1)\in\mbox{DNA}(\alpha) \mbox{ (see Mikosch \textit{et al.} \cite{MGKA:1995})},
\end{eqnarray*}
where DNA$(\alpha)$ denotes the set of random variables in the domain
of normal attraction of
a symmetric $\alpha$-stable random variable.
Letting $\bm{\beta}=(\phi_1,\ldots,\phi_p,\theta_1,\ldots,\theta
_q)$, define
\begin{eqnarray*}
&&C = \bigl\{ \bm{\beta}\in\R^{p+q} \dvtx  \phi_p,
\theta_q \neq0, \phi(z) \mbox{ and}
\\
&&\hphantom{C = \bigl\{}\theta(z)\mbox{ have no common zeros, }\phi(z)\theta(z)\neq0 \mbox{ for }
\llvert z\rrvert \leq1 \bigr\},
\end{eqnarray*}
where $\phi(z) = 1+\phi_1 z + \cdots+ \phi_p z^p$, and $\theta(z)
= 1
+ \theta_1 z + \cdots+ \theta_q z^q$. Let $g(\omega;\bm{\beta})$ be
\[
g(\omega;\bm{\beta}) = \biggl\llvert \frac{1+\theta_1 \exp(\mathrm{i}\omega) +
\cdots+
\theta_q \exp(\mathrm{i} q\omega)}{1+\phi_1 \exp(\mathrm{i}\omega) + \cdots+ \phi
_p \exp
(\mathrm{i} p\omega)} \biggr\rrvert
^2.
\]
They defined the Whittle estimator of $\bm{\beta}$ by
\[
\hat{\bm{\beta}}_n \equiv\arg\min_{\bm{\beta}\in C}\int
^{\uppi
}_{-\uppi} \frac{\tilde
{I}_{n,X}(\omega)}{g(\omega;\bm{\beta})}\,\mathrm{d}\omega.
\]
Then Mikosch \textit{et al.} \cite{MGKA:1995} showed that the estimator $\hat{\bm{\beta
}}_n$ is consistent to the true parameter $\bm{\beta}_0\in C$.

In many cases, however, we know neither the true stochastic structure
of the process
nor the true pivotal unknown quantities.
In such cases, we can apply the empirical likelihood approach to the data,
without assuming that the data come from a known family of stochastic models.
The empirical likelihood approach was introduced as a non-parametric
method of inference based on a data-driven likelihood ratio function
in the i.i.d. case (e.g., Owen \cite{O:1988}).
For dependent data, Monti \cite{M:1997} applied the empirical likelihood
approach to a stationary linear process with the finite second moment.
She used
\[
\bmm{m}(\lambda_t;\bm{\theta}) = \frac{\partial}{\partial\bm
{\theta
}} \biggl\{ \log g(
\lambda_t;\bm{\theta}) + \frac{I_{n, X}(\lambda_t)}{g(\lambda_t;\bm{\theta})} \biggr\},\qquad t=1,\ldots, n
\]
as an \textit{estimating function}.
This can be understood as a discretized derivative of the Whittle likelihood
\[
\int^{\uppi}_{-\uppi} \biggl\{ \log g(\omega;\bm{\theta})
+ \frac{I_{n, X}(\omega)}{g(\omega;\bm{\theta})} \biggr\}\,\mathrm{d}\omega.
\]
Here $g(\omega;\bm{\theta})$ and $I_{n, X}(\omega)$ are, respectively,
the usual spectral density of a stationary process and the periodogram.
Using this estimating function, the empirical likelihood ratio function
is defined as
%
\begin{eqnarray}
R(\bm{\theta}) = \max_{w_1,\dots, w_n} \Biggl\{ \prod
_{t=1}^n n w_t; \sum
_{t=1}^n w_t \bmm{m}(
\lambda_t; \bm {\theta}) = \bm{0}, \sum_{t=1}^n
w_t = 1, 0 \leq w_1, w_2, \dots,
w_n \leq1 \Biggr\}.\qquad\label{eq:2.2}
\end{eqnarray}
Under the circular assumption,
It is shown that the quantity $-2\log R(\bm{\theta})$ converges in
distribution
to chi-square random variable with degree of freedom $q$
under $H$: $\bm{\theta}=\bm{\theta}_0$ (the pivotal true value of
$\bm
{\theta}$) $\in\Theta\subset\R^q$.
Ogata and Taniguchi \cite{OT:2010} developed the empirical likelihood approach
to multivariate non-Gaussian stationary processes
without the circular assumption.
For a vector-valued process $\{\bmm{X}(t); t \in\Z\}$,
\[
\bmm{X}(t) = \sum_{j=0}^{\infty} \bmm{G}(j)
\bmm{e}(t-j),\qquad E\bigl[\bmm {e}(t)\bmm{e}(l)'\bigr] = \delta(t, l)
\bm{\Sigma},
\]
they introduced the disparity measure
\[
D(\bmm{f}_{\bm{\theta}};\bmm{g}) = \int^{\uppi}_{-\uppi}
\bigl[\log \det\bmm {f}(\omega;\bm{\theta}) + \tr\bigl\{\bmm{f}(\omega;\bm{
\theta})^{-1}\bmm {g}(\omega)\bigr\} \bigr]\,\mathrm{d}\omega
\]
on
\[
\mathcal{P} = \Biggl\{ \bmm{f}(\omega;\bm{\theta}) |\bmm {f}(\omega;\bm {
\theta}) = \Biggl\{\sum^{\infty}_{j=0}\bmm{G}(j;
\bm{\theta})\exp (\mathrm{i}j\omega) \Biggr\} \bm{\Sigma} \Biggl\{\sum
^{\infty}_{j=0}\bmm{G}(j;\bm{\theta})\exp(\mathrm{i}j\omega )
\Biggr\} ^\ast, \bm{\theta}\in\Theta\subset\R^q \Biggr\},
\]
where $\bmm{g}(\omega)$ is the usual spectral density matrix of the
$s$-dimensional stationary linear process.
If the innovation variance of the process is independent of unknown
parameter $\bm{\theta}$, we call $\bm{\theta}$ ``\textit{innovation
free}''. Then,
the first integration of the disparity measure is
independent of $\bm{\theta}$
(e.g., Hannan \cite{hannan2009multiple}, page 162).
Therefore if $\bm{\theta}$ is innovation-free, the derivative of this
measure is
\[
\frac{\partial}{\partial\bm{\theta}}D(\bmm{f}_{\bm{\theta
}};\bmm{g}) = \frac{\partial}{\partial\bm{\theta}}\int
^{\uppi}_{-\uppi}\tr\bigl\{\bmm {f}(\omega;\bm {
\theta})^{-1}\bmm{g}(\omega)\bigr\}\,\mathrm{d}\omega.
\]
They introduced the pivotal true value $\bm{\theta}_0$ defined by
%
\begin{eqnarray}\label{eq:2.5}
\frac{\partial}{\partial\bm{\theta}}\int^{\uppi}_{-\uppi}\tr\bigl\{ \bmm
{f}(\omega;\bm {\theta})^{-1}\bmm{g}(\omega) \bigr\}\,\mathrm{d}\omega
\Big|_{\bm{\theta}=\bm
{\theta
}_0} = \bm{0}.
\end{eqnarray}
In this case, the estimating function is naturally set to be
\[
\bmm{m}(\lambda_t;\bm{\theta}) = \frac{\partial}{\partial\bm
{\theta}}\tr \bigl\{
\bmm{f}(\lambda_t;\bm{\theta})^{-1} \bmm{I}_{n, X}(
\lambda_t) \bigr\},\qquad t=1,\ldots,n,
\]
where $\bmm{I}_{n, X}(\omega)$ is the usual periodogram matrix.
Under mild conditions on the fourth order cumulant of the process, they
showed that $-2\log R(\bm{\theta})$ converges in law to a sum of gamma
distributed random variables under $H$: $\bm{\theta}=\bm{\theta}_0$.

The approach has been discussed for stationary processes with the
``finite second moments''.
In this paper, we consider a linear process $\{X(t); t\in\Z\}$
generated by \eqref{eq:1.1}
with $\{ Z(t); t\in\Z\}$, a sequence of i.i.d. symmetric $\alpha
$-stable random variables with scale $\sigma>0$,
and the characteristic function of $Z(1)$ is given as
\[
\E\exp\bigl\{\mathrm{i} Z(1) \xi\bigr\} = \exp\bigl\{-\sigma\llvert \xi\rrvert
^\alpha\bigr\},\qquad\xi\in\R.
\]
Generally, we can define the stable process for $\alpha\in(0,2]$.
However, we assume that $\alpha\in[1,2)$ to guarantee probability
convergence of
important terms which will appear in proofs of theorems in this paper.
This restriction is not quite strict, since the process \eqref{eq:1.1}
with $\alpha\in[1,2)$ still does not have the
finite second moment.
To guarantee the a.s. absolute convergence of \eqref{eq:1.1}, we make
the following assumption.

\begin{ass}\label{ass:2.2}
For some $\delta$ satisfying $0<\delta<1$,
\[
\sum_{j=0}^\infty\llvert j\rrvert \llvert
\psi_j\rrvert ^\delta<\infty.
\]
\end{ass}

Under this assumption, the series \eqref{eq:1.1} converges almost
surely. This is an easy consequence of the three-series theorem (c.f.
Petrov \cite{P:1975}).
Furthermore, the process \eqref{eq:1.1} has the normalized power
transfer function
\[
\tilde{g}(\omega) =\frac{1}{\psi^2}{ \Biggl\llvert \sum
^{\infty
}_{j=0}\psi_j \exp(\mathrm{i}j\omega) \Biggr
\rrvert ^2},\qquad\psi^2 = \sum^{\infty}_{j=0}
\psi_j^2.
\]
From the property of stable random variables,
\[
X(t) =_d \Biggl\{\sum^\infty_{j=0}
\llvert \psi_j\rrvert ^\alpha \Biggr\} ^{1/\alpha} Z(1),
\]
which implies that this process does not have the finite second moment
when $\alpha<2$,
so we cannot use the method of moments.
The empirical likelihood approach is still useful when we deal with the
stable process. Hereafter, we define a
pivotal true value $\bm{\theta}_0$ of the process \eqref{eq:1.1} as the
solution of
%
\begin{eqnarray}\label{eq:11}
\frac{\partial}{\partial\bm{\theta}}\int^{\uppi}_{-\uppi}
\frac
{\tilde
{g}(\omega
)}{f(\omega;\bm{\theta})}\,\mathrm{d}\omega \Big|_{\bm{\theta}=\bm{\theta
}_0} = \bm {0},
\end{eqnarray}
where $\bm{\theta}=(\theta_1,\ldots,\theta_p)'\in\Theta\subset
\R^p$.
Note that the score function does not necessarily coincide with the
true normalized power transfer function $\tilde{g}(\omega)$,
and we can choose various important quantities $\bm{\theta}_0$ by
choosing the form of $f(\omega;\bm{\theta})$.
For example, for fixed $l\in\N$, set
\[
f(\omega;\theta) = \bigl\llvert 1-\theta\exp(\mathrm{i}l\omega)\bigr\rrvert
^{-2}.
\]
Solving (\ref{eq:11}), we have
\[
{\theta}_0 = \frac{\sum^{\infty}_{j=0}\psi_j\psi_{j+l}}{\sum^{\infty
}_{j=0}\psi_j^2} \equiv\rho(l) \qquad(\mbox{say}).
\]
On the other hand, a sample autocorrelation function
\[
\hat{\rho}(l) \equiv\frac{\sum^{n-l}_{t=1}X(t)X(t+l)}{\sum^{n}_{t=1}X(t)^2},\qquad l\in\N
\]
for the stable process (\ref{eq:1.1}) is weakly consistent to the
autocorrelation function of the process in the stable case; namely,
for fixed $l$,
$\prolim_{n\to\infty}\hat{\rho}(l)=\rho(l)$
(e.g., Davis and Resnick \cite{DR1986}).

By these motivation, we consider the empirical likelihood ratio
function \eqref{eq:2.2}
with
\[
{\bmm{m}}(\lambda_t;\bm{\theta}) = \frac{\partial}{\partial\bm
{\theta
}}
\frac{\tilde{I}_{n,X}(\lambda_t)}{f(\lambda_t;\bm{\theta})},\qquad \lambda_t = \frac{2\uppi t}{n}, t=1,\ldots,n.
\]
Hereafter, we make the following assumptions on $f(\omega;\bm{\theta})$.

\begin{ass}\label{ass:2.1}\mbox{}
\begin{enumerate}[(iii)]
\item[(i)]$\Theta$ is a compact subset of $\R^q$ and
$f(\omega;\bm{\theta})$ has an parametrized representation
as an element of $\mathcal{P}$, where $\mathcal{P}$ is defined by
\[
\mathcal{P} = \Biggl\{ f(\omega;\bm{\theta}) |{f}(\omega;\bm {\theta}) = \Biggl| \sum
^{\infty}_{j=0}\eta_j(\bm{\theta})
\exp(\mathrm{i}j\omega) \Biggr|^2, \bm{\theta}\in\Theta\subset\R^q \Biggr
\}.
\]
\item[(ii)] For any $\bm{\theta} \in\Int\Theta$,
$f(\omega;\bm{\theta})$ is continuously twice differentiable with
respect to $\bm{\theta}$.
\item[(iii)] There exists an unique $\bm{\theta}_0\in\Theta$ satisfying
\eqref{eq:11}.
\end{enumerate}
\end{ass}

\section{Main results}\label{sec3}
In this section, we introduce the limit distribution of the empirical
likelihood statistic
for the scalar stable process (\ref{eq:1.1}).
Our main purpose is to make an accurate confidence region of
${\bm{\theta}}_0$ based on the empirical likelihood approach.
Because of the properties of stable random variables,
it is difficult to use the method of moments.
To overcome this problem, we frequently make use of the self-normalized
periodogram defined in Section~\ref{sec2}.
Kl\"{u}ppelberg and Mikosch \cite{KM:1994} or Mikosch \textit{et al.} \cite{MGKA:1995}
introduced the self-normalized periodogram,
and Kl\"{u}ppelberg and Mikosch \cite{KM:1996} showed some limit
theorems of integrated self-normalized periodogram.
Under the settings in Section~\ref{sec2},
we derive the asymptotic distribution of the empirical likelihood ratio
statistic,
and construct a confidence region for $\bm{\theta}_0$.

We impose an assumption to describe the asymptotics of the empirical
likelihood ratio statistic.

\begin{ass} \label{ass:3.1}
For some $\mu\in(0,\alpha)$ and all $k=1,\ldots,q$,
\[
\sum^{\infty}_{t=1} \biggl\llvert \int
^{\uppi}_{-\uppi}\frac{\partial
}{\partial
\theta_k}\frac{\tilde{g}(\omega)}{{f}(\omega;\bm{\theta})}
\Big| _{\bm
{\theta}=\bm{\theta}_0} \cos( t\omega)\,\mathrm{d}\omega \biggr|^{\mu}<\infty.
\]
\end{ass}

Assumption~\ref{ass:3.1} is used for Proposition~3.5 of Kl\"{u}ppelberg and Mikosch \cite{KM:1996}.
It is easy to see that stable AR($p$) processes
satisfying Assumption~\ref{ass:2.1}
satisfy this assumption.

In order to control the rate of convergence of the empirical likelihood
ratio statistic, we introduce the normalizing sequence
\[
x_n = \biggl(\frac{n}{\log n} \biggr)^{1/\alpha},\qquad n=2,3,
\ldots.
\]

The next theorem gives the asymptotics of $R(\bm{\theta}_0)$.
The proof will be given in Section~\ref{sec6}.

\begin{theorem}\label{theorem:1}
Suppose that $\alpha\in[1,2)$, and Assumptions \ref{ass:2.2}, \ref
{ass:2.1} and \ref{ass:3.1} hold.
Then,
%
\begin{eqnarray}\label{eq:3.1}
-\frac{2 x_n^2}{n}\log R(\bm{\theta}_0) \dlim
\bmm{V}'\bmm{W}^{-1}\bmm{V} \qquad\mbox{under } H:
\bm{\theta}=\bm{\theta}_0,
\end{eqnarray}
where
$\bmm{V}$ and $\bmm{W}$ are $q\times1$ random vector and $q\times q$
constant matrix, respectively, whose $j$th and $(k,l)$-elements are
expressed as
\begin{eqnarray*}
V_j &=& \frac{1}{\uppi}\sum^{\infty}_{t=1}
\frac{S_t}{S_0} \biggl\{ \int^{\uppi}_{-\uppi}
\frac{\partial f(\omega;\bm{\theta})^{-1}}{\partial\theta_j}\Big |_{\bm{\theta}=\bm{\theta}_0}\tilde{g}(\omega)\cos( t\omega)\,\mathrm{d}\omega \biggr
\},
\\
W_{k l} &=& \frac{1}{2\uppi}\int^{\uppi}_{-\uppi}
\frac{\partial
f(\omega;\bm
{\theta
})^{-1}}{\partial\theta_k} \frac{\partial f(\omega;\bm{\theta})^{-1}}{\partial\theta_l} \Big|_{\bm{\theta}=\bm{\theta}_0} 2 \tilde{g}(
\omega)^2\,\mathrm{d}\omega
\end{eqnarray*}
with independent random variables $S_0$, $S_1$, $S_2,\ldots$; $S_0$ is
a positive $\alpha/2$-stable random variable and $\{S_j; j=1,2,\ldots
\}$
is a sequence of symmetric $\alpha$-stable random variables.
\end{theorem}

\begin{rem}
The limit distribution (\ref{eq:3.1}) depends on the characteristic
exponent $\alpha$
and unknown normalized power transfer function $\tilde{g}(\omega)$.
We can construct appropriate consistent estimators of them.
It is shown that Hill's estimator
\[
\hat{\alpha}_{\mathrm{Hill}} = \Biggl\{\frac{1}{k}\sum
^{k}_{t=1}\log \frac
{\llvert X\rrvert _{(t)}}{\llvert X\rrvert _{(k+1)}} \Biggr
\}^{-1}
\]
is a consistent estimator of $\alpha$, where $\llvert X\rrvert _{(1)}>\cdots
>\llvert X\rrvert _{(n)}$ is the order statistic of $\llvert X(1)\rrvert,\ldots,\allowbreak\llvert X(n)\rrvert $
and $k=k(n)$ is an integer satisfying some conditions
(e.g., Resnick and St\v{a}ric\v{a} \cite{RS:1996} and \cite{RS:1998}).
Next, it is known that the smoothed self-normalized periodogram by an
appropriate weighting function $W_n(\cdot)$ is
weekly consistent to the normalized power transfer function. That is,
\[
\tilde{J}_{n,X}(\omega) = \sum_{\llvert k\rrvert \leq m}W_n(k)
\tilde {I}_{n,X}(\lambda _k)\plim \tilde{g}(\omega),\qquad
\lambda_k = \omega+ \frac{k}{n}, \llvert k\rrvert \leq m
\]
for any $\omega\in[-\uppi,\uppi]$ (Kl\"{u}ppelberg and Mikosch
\cite{KM:1993}, Theorem~4.1), where the integer $m=m(n)$
satisfies $m\to\infty$ and $m/n\to0$ as $n\to\infty$.
One possible choice of the weighting function $W_n(\cdot)$ and $m=m(n)$
are $W_n(k)=(2m+1)^{-1}$ and $m=[\sqrt{n}]$ ($[x]$ denotes the integer
part of $x$).
We use this weighting function in the section of numerical studies.
Then, by Slutsky's lemma and continuous mapping theorem, we obtain
consistent estimator
$\hat{\bmm{W}}$ of $\bmm{W}$.
So if we choose a proper threshold value $\gamma_p$, which is
$p$-percentile corresponding to $\bmm{V}'\bmm{W}\bmm{V}$,
$C_{\alpha, p}$ below is an approximate $p/100$ level confidence region
of $\bm{\theta}_0$.
%
\begin{eqnarray}
C_{\alpha, p} = \biggl\{ \bm{\theta}\in\Theta; -\frac{2
x_n^2}{n}\log R(
\bm{\theta})<\gamma_p \biggr\}.\label{eq:3.2}
\end{eqnarray}
\end{rem}

\section{Vector \texorpdfstring{$\alpha$}{$alpha$}-stable processes}\label{sec4}
So far we focused on the scalar case for clarity.
In this section, we extend the empirical likelihood analysis to
the case of vector $\alpha$-stable processes.
Consider a $d$-dimensional vector-valued linear process $\{ \bmm{X}(t);
t \in\Z\}$ generated by
%
\begin{equation}
\label{eq:4.1} \bmm{X}(t) = \sum_{j=0}^{\infty}
\Psi(j) \bmm{Z}(t-j),
\end{equation}
where $\Psi(0)$ is the identity matrix and $\{\Psi(j); j\in\N\}$ is a
sequence of $d \times d$ real matrices, and
$\{ \bmm{Z}(t); t \in\Z\}$ is an independently and identically
distributed sequence of
symmetric $\alpha$-stable random vectors whose elements
are also independent.

Now, we set down the following assumptions for the general result.
Almost all assumptions are similar to those of the 1-dimensional stable
processes.

\begin{ass}\label{ass:4.0}
For some $\delta$ satisfying $0<\delta<1$ and all $k, l = 1,\dots, d$,
%
\begin{equation}
\label{eq:4.2} \sum_{j=0}^{\infty} j\bigl
\llvert \Psi(j)_{kl}\bigr\rrvert ^{\delta} < \infty.
\end{equation}
\end{ass}

The sample autocovariance and the periodogram matrices are defined as
\begin{eqnarray*}
\hat{\Gamma}_{n, X}(h) &= &n^{-2/\alpha} \sum
_{t=1}^{n-\llvert h\rrvert }\bmm{X}(t) \bmm{X}(t+h)',
\\
\bmm{I}_{n, X}(\omega) &=& d_{n, X}(\omega) d_{n, X}(
\omega)^*,\qquad d_{n, X}(\omega)= n^{-1/\alpha} \sum
_{t=1}^n \bmm{X}(t) \exp (\mathrm{i}\omega t),
\end{eqnarray*}
respectively.
We define the true power transfer function $\bmm{g}(\omega)$ by
\[
\bmm{g}(\omega) = \Psi(\omega)\Psi(\omega)^*,
\]
where $\Psi(\omega) = \sum_{j=0}^{\infty} \Psi(j) \exp(\mathrm{i}j\omega)$.
Similarly as in the previous section,
we use the empirical likelihood ratio with the estimating function
\[
{\bmm{m}}(\lambda_t; \theta) = \frac{\partial}{\partial\bm{\theta}} \tr\bigl\{
\bmm{f}(\lambda_t; \bm{\theta})^{-1} \bmm{I}_{n, X}(
\lambda_t) \bigr\},
\]
where $\bmm{f}(\lambda_t; \bm{\theta})$ satisfies the following assumptions.

\begin{ass}\label{ass:4.1}\mbox{}
\begin{enumerate}[(iii)]
\item[(i)]$\Theta$ is a compact subset of $\R^q$ and
$\bmm{f}(\omega;\bm{\theta})$ has an parametrized representation
as an element of $\mathcal{P}$, where $\mathcal{P}$ is defined by
\[
\hspace*{-8pt}\mathcal{P}= \Biggl\{\bmm{f}(\omega; \bm{\theta}) | \bmm{f}(\omega; \bm {
\theta}) = \Biggl(\sum_{j=0}^{\infty}\Xi(j; \bm{
\theta}) \exp(\mathrm{i}j\omega ) \Biggr) \Biggl(\sum_{j=0}^{\infty}
\Xi(j; \bm{\theta}) \exp(\mathrm{i}j \omega) \Biggr)^*, \bm{\theta} \in\Theta\subset
\R^q \Biggr\}.
\]
\item[(ii)] For any $\bm{\theta} \in\Int\Theta$,
$f(\omega;\bm{\theta})$ is continuously twice differentiable with
respect to $\bm{\theta}$.
\item[(iii)] There exists an unique $\bm{\theta}_0\in\Theta$ satisfying
\eqref{eq:11}.
\end{enumerate}
\end{ass}

Assumption~\ref{ass:4.2} below guarantees the convergence of
the functional of periodogram by inequality of an application
of Theorem~3.1 in Rosinski and Woyczynski \cite{rosinski1987multilinear}.

\begin{ass} \label{ass:4.2}
For some $\mu\in(0,\alpha)$ and all $k=1,\ldots,q$,
\[
\sum^{\infty}_{t=1} \biggl\llVert \int
^{\uppi}_{-\uppi}\frac{\partial
}{\partial\theta_k} \Psi(
\omega)^{*} \bmm{f}(\omega; \bm{\theta}) \Psi(\omega) \exp (\mathrm{i} t \omega)\,\mathrm{d}\omega \biggr\rrVert ^{\mu}_{E} < \infty.
\]
\end{ass}

\begin{theorem}\label{theorem:4.1}
Suppose that $\alpha\in[1,2)$, and Assumptions \ref{ass:4.0}--\ref
{ass:4.2} hold
for the process \eqref{eq:4.1}.
If
\[
\frac{\partial}{\partial\bm{\theta}} \int_{-\uppi}^{\uppi} \Psi(\omega)^*
\bmm{f}(\omega;\bm{\theta })^{-1} \Psi (\omega)\,\mathrm{d}\omega \Big
\vert_{\bm{\theta} = \bm{\theta_0}} = \bm{0},
\]
then
\[
-2 \frac{x_n^2}{n} \log R(\bm{\theta}_0) \dlim
\bmm{V}' \bmm{W}^{-1} \bmm{V}\qquad \mbox{under $H$: $\bm{
\theta}=\bm {\theta}_0$,}
\]
where
\[
\bmm{V} = \frac{1}{2\uppi} \sum_{i, j =1}^d
\sum_{h=1}^{\infty} \frac{S(h)_{ij}}{S_{\alpha/2}}
\pmatrix{\displaystyle \int_{-\uppi}^{\uppi}
\bigl(B_1(\omega) + \overline{B_1(\omega)}
\bigr)_{ij}\,\mathrm{d}\omega
\cr
\displaystyle\int_{-\uppi}^{\uppi} \bigl(B_2(\omega) +
\overline{B_2(\omega)}\bigr)_{ij}\,\mathrm{d}\omega
\cr
\vdots
\cr
\displaystyle\int_{-\uppi}^{\uppi} \bigl(B_q(\omega) +
\overline{B_q(\omega)}\bigr)_{ij}\,\mathrm{d}\omega, }
\]
with $S(h)_{ij}$ a matrix whose all elements are stable with index
$\alpha$,
$S_{\alpha/2}$ a random variable with index $\alpha/2$ and
\[
B_k(\omega) = \Psi(\omega)^* \frac{\partial}{\partial\theta_k} \bm {f}(\omega;
\bm{\theta})^{-1} \Big\vert_{\bm{\theta} = \bm{\theta_0}} \Psi(\omega)\exp(\mathrm{i} h
\omega),\qquad k=1, \dots, q,
\]
and the $(a, b)$-component of $\bmm{W}$ can be expressed as
%
\begin{eqnarray*}
&&W_{ab}= \frac{1}{2\uppi d^2} \int^{\uppi}_{-\uppi}
\biggl( \tr \biggl[ \bmm{g}(\omega) \frac{\partial\bmm{f}(\omega; \bm{\theta})^{-1}}{\partial\theta
_a}\Big|_{\bm{\theta} = \bm{\theta_0}}
\bmm{g}(\omega) \frac{\partial\bmm{f}(\omega; \bm{\theta})^{-1}}{\partial\theta
_b}\Big|_{\bm{\theta} = \bm{\theta_0}} \biggr]
\\
&&\hphantom{W_{ab}= \frac{1}{2\uppi d^2} \int^{\uppi}_{-\uppi}
\biggl( }{}+ \tr \biggl[ \bmm{g}(\omega) \frac{\partial\bmm{f}(\omega; \bm{\theta})^{-1}}{\partial\theta
_a}\Big|_{\bm{\theta} = \bm{\theta_0}} \biggr]
\tr \biggl[ \bmm{g}(\omega) \frac{\partial\bmm{f}(\omega; \bm{\theta})^{-1}}{\partial\theta
_b}\Big|_{\bm{\theta} = \bm{\theta_0}} \biggr]
\biggr)\,\mathrm{d}\omega.
\end{eqnarray*}
\end{theorem}

\begin{pf}
The proof of Theorem~\ref{theorem:4.1} is given in the supplemental
article (Akashi \textit{et al.} \cite{akashi}),
since it is more technical.
\end{pf}

\begin{rem}
This extension is not straightforward, and
contains some novel aspects.
We take up an appealing example for Theorem~\ref{theorem:4.1}.
Consider whether the wave structures
of the spectra between all components are
``close'' to each other or not.
For simplicity, we formulate this idea in 2-dimensional case
and assume the true power transfer function $\bmm{g}(\omega)$ is
\[
\bmm{g}(\omega) = \frac{1}{2\uppi} \sum_{k=-\infty}^{\infty}
\tilde {\bmm{R}}(k) \exp(-\mathrm{i} k\omega),
\]
where $\tilde{\bmm{R}}(k)$, a symmetric matrix, denotes the $k$th
autocorrelation function.
Then the null hypothesis can be written as
\[
H\dvt  \tilde{\bmm{R}}(k) = \theta_0 \tilde{\bmm{R}}(j)\quad \mbox{or}\quad
\tilde{\bmm{R}}(k) = \theta_0 \tilde{\bmm{R}}(j)'
\qquad\mbox{for some $k$ and $j$}.
\]
To test this hypothesis, we set the estimating function ${\bmm
{m}}(\lambda_t; \bm{\theta})$
with an inverse correlation function $\bmm{f}(\lambda_t; \bm{\theta})^{-1}$,
which was first introduced in Cleveland \cite{cleveland1972inverse}, and deeply discussed by
Bhansali~\cite{bhansali1980autoregressive}.
Let
\[
\bmm{f}(\omega; \theta)^{-1} = \bigl(\exp(k \omega) + \exp(-k \omega)
\bigr) %
\pmatrix{ \theta& 0
\cr
0 & \theta } %
+ \bigl(\exp(j
\omega) + \exp(-j \omega)\bigr) %
\pmatrix{ \frac{1}{2}
\theta^2 & 0
\cr
0 & \frac{1}{2}\theta^2 }
.
\]
Then under the hypothesis, we have
\[
\frac{\partial}{\partial\theta} \int_{-\uppi}^{\uppi} \Psi(\omega)^*
\bmm{f}(\omega;\bm{\theta })^{-1} \Psi (\omega)\,\mathrm{d}\omega \Big
\vert_{\theta= \theta_0} = \bm{0},
\]
which satisfies the assumption in Theorem~\ref{theorem:4.1}.
\end{rem}

\section{Numerical studies}\label{sec5}
In this section, we carry out some simulation studies for Theorems \ref{theorem:1}
and \ref{theorem:4.1}.
Suppose that the observations $X(1),\ldots, X(n)$ are generated from
the following scalar-valued stable MA($100$) model:
%
\begin{eqnarray}
X(t) = \sum^{100}_{j=0}\psi_j
Z(t-j), \label{5.1}
\end{eqnarray}
where $\{Z(t); t\in\Z\}$ is a sequence of i.i.d. s$\alpha$s random
variables with scale $\sigma=1$ and coefficients $\{\psi_j; j\in\N\}$
are defined as
\[
\psi_j= \cases{ 1 &\quad $(j=0)$,
\cr
b^j/j &\quad $(1\leq j\leq100)$,
\cr
0 &\quad (otherwise).} %
\]
Since this process can not be expressed as AR or ARMA models with
finite dimension,
it is suitable to apply the empirical likelihood approach to estimate
pivotal unknown quantities.
We first discuss the estimation of the autocorrelation with lag $2$
%
\begin{eqnarray}
\rho(2) = \prolim_{n\to\infty}\frac{\sum^{n-2}_{t=1}X(t)X(t+2)}{\sum^{n}_{t=1}X(t)^2}.\label{5.3}
\end{eqnarray}
It is seen that the normalized power transfer function of the process
(\ref{5.1}) is given by
\[
\tilde{g}(\omega) = \frac{ \llvert  \sum^{100}_{j=0}\psi_j\exp
(\mathrm{i}j\omega)
 \rrvert ^2}{\sum^{100}_{j=0}\psi_j^2}.
\]
If we set the score function as $f(\omega;\theta) = \llvert 1-\theta\exp
(2\mathrm{i}\omega)\rrvert ^{-2}$,
we obtain
\[
\theta_0=\frac{\sum^{100}_{j=0}\psi_j\psi_{j+2}}{\sum^{100}_{j=0}\psi_j^2}.
\]
On the other hand, from Davis and Resnick \cite{DR1986}, the right-hand side
limit of (\ref{5.3}) exists, and is equal to this $\theta_0$.
So it is natural that we define the estimating function $m(\lambda
_t;\theta)$ by this $f(\omega;\theta)$ to estimate $\rho(2)$.
The autocorrelation can also be estimated by sample autocorrelation
(SAC) method.
From Theorem~12.5.1 of Brockwell and Davis \cite{BD:1991}, for fixed $l\in\N$,
\[
x_n\bigl\{ \hat\rho(l) - \rho(l) \bigr\} \dlim\frac{\widetilde S_1}{\widetilde S_0}
\Biggl\{\sum^{\infty
}_{j=1}\bigl\llvert
\rho(l+j) + \rho(l-j) - 2\rho(j)\rho(l) \bigr\rrvert ^\alpha \Biggr\}
^{1/\alpha},
\]
where $\hat\rho(l) = {\sum^{n-l}_{t=1}X(t)X(t+l)}/{\sum^{n}_{t=1}X(t)^2}$, $\widetilde S_0$ and $\widetilde S_1$ are $\alpha
/2$ and $\alpha$-stable random variables, respectively.
Under this setting, we construct confidence intervals of $\theta_0 =
\rho(2)$ by calculating $R(\theta)$ at numerous point over $(-1,1)$,
and compare confidence intervals constructed by the empirical
likelihood method with the SAC method.

The results of our simulations are as follows.
First, we generate $300$ samples from (\ref{5.1}). Note that in this
case, the characteristic exponent $\alpha=1.5$ is known.
Then using the weighting function $W_n$ which is mentioned in Section~\ref{sec3},
we calculate the consistent estimator $\tilde{J}_{n,X}(\omega)$ of
$\tilde{g}(\omega)$
and construct an approximate $90$\% confidence interval of ${\theta}_0$
defined as \eqref{eq:3.2}.
We also use the Monte Carlo simulation to calculate $\gamma_{90}$ which
is 90 percentile of $\bmm{V}'\bmm{W}\bmm{V}$
for $10^5$ times.
Table~\ref{table:1} shows the values of $\theta_0$ and confidence
intervals by
the empirical likelihood method and the sample autocorrelation method
for $b=0.5$ (case 1) and $0.9$ (case 2).
By this simulation, it is shown that the length of intervals
obtained by the empirical likelihood method is seems to be shorter than
that by the sample autocorrelation
method.

\begin{table}[b]
\tabcolsep=0pt
\caption{90\% confidence intervals (and length) for the autocorrelation
with lag $2$. Sample size is $300$ and $\alpha=1.5$}\label{table:1}
\begin{tabular*}{\textwidth}{@{\extracolsep{\fill}}llllllll@{}}
\hline
& $\theta_0\approx$ & \multicolumn{3}{l}{E.L.} & \multicolumn{3}{l}{SAC} \\
\hline
Case 1 & 0.1168 & $-$0.0761 & 0.1930 & (0.2691) & $-$0.0676 & 0.2481 & (0.3157) \\
Case 2 & 0.3603 & \hphantom{$-$}0.1320 & 0.4765 & (0.3445) & \hphantom{$-$}0.1388 & 0.5304 & (0.3916) \\
\hline
\end{tabular*}
\end{table}

Next, we fix $b = 0.5$ and $n=300$, and construct confidence intervals
for cases of $\alpha=1.0$ (Cauchy), $1.5$ and $1.9$ (near Gaussian).
The larger $\alpha$ becomes, the better performance both
methods show (see Table~\ref{table:2}). In particular, the empirical
likelihood method provides better
inferences than those by the SAC method when $\alpha$ is nearly $1$.

\begin{table}[t]
\tabcolsep=0pt
\caption{90\% confidence intervals (and length) for the autocorrelation
with lag $2$. Sample size is $300$, $b= 0.5$ and $\theta_0\approx
0.1168$}\label{table:2}
\begin{tabular*}{\textwidth}{@{\extracolsep{\fill}}llllllll@{}}
\hline
&$\alpha$ & \multicolumn{3}{l}{E.L.} & \multicolumn{3}{l}{SAC} \\
\hline
Case 3 & 1.0 & $-$0.1583 & 0.3335 & (0.4918) & $-$0.1342 & 0.3891 & (0.5233) \\
Case 4 & 1.5 & $-$0.0761 & 0.1930 & (0.2691) & $-$0.0676 & 0.2481 & (0.3157) \\
Case 5 & 1.9 & $-$0.0465 & 0.1329 & (0.1794) & $-$0.0450 & 0.1365 & (0.1815) \\ \hline
\end{tabular*}
\end{table}

Moreover, we investigate the length of intervals when $b=0.5$ and
$\alpha=1.5$ for small samples.
Table~\ref{table:3} shows the result for $n=50$ and $100$.
Even though sample size is small,
the empirical likelihood method also works well.

\begin{table}[b]
\tabcolsep=0pt
\caption{90\% confidence intervals (and length) for the autocorrelation
with lag $2$. $b = 0.5$, $\alpha=1.5$
and $\theta_0\approx0.1168$}\label{table:3}
\begin{tabular*}{\textwidth}{@{\extracolsep{\fill}}llllllll@{}}
\hline
& $n$ & \multicolumn{3}{l}{E.L.} & \multicolumn{3}{l}{SAC} \\ \hline
Case 6 &\hphantom{1}50 & $-$0.2397 & 0.4313 & (0.6710) & $-$0.2477 & 0.5629 & (0.8106) \\
Case 7 &100 & $-$0.3125 & 0.2228 & (0.5353) & $-$0.3476 & 0.2218 & (0.5694) \\ \hline
\end{tabular*}
\end{table}

Also, we give an example for multivariate case.
Suppose that the observations $\bmm{X}(1), \ldots, \bmm{X}(n)$
are generated from the 2-dimensional VMA(100) model with
innovations $\{\bmm{Z}(t); t\in\Z\}$
whose marginal distributions are i.i.d. s$\alpha$s with scale 1,
and the coefficient matrices $A(j)$, $j=1,\dots, 100$ are assumed to be
\[
A(j)= %
\pmatrix{ 0.7^j & j^{-2} b^j
\cr
0 & 0.5^j } %
.
\]
To this model, we use the following score function $\bmm{f}(\omega;
\theta)$ defined by
\[
\bmm{f}(\omega; \theta) = \bigl(I-B_{\theta} \exp(\mathrm{i} \omega)
\bigr)^{-1} {\bigl(I-B_{\theta} \exp(\mathrm{i} \omega)
\bigr)^{-1}}^*, \qquad\mbox{where } B_{\theta}= %
\pmatrix{ 0.5
& \theta
\cr
0.4 & 0.2 } %
.
\]
In this case, the asymptotic distribution of $-2({x_n^2}/{n}) \log
R(\theta_0)$
can be simply represented by $ ({S_1}/{S_0}
)^2({V^2}/{W})$, where
$S_0$ and $S_1$ are the same as in Theorem~\ref{theorem:1}, $W$ is the
same as
in Theorem~\ref{theorem:4.1} and
\[
V = \frac{1}{\uppi} \Biggl[ \Biggl\llvert \int_{-\uppi}^{\uppi}
F_{12}(\omega) \cos( \omega)\,\mathrm{d}\omega + \sum
_{t=1}^{\infty} \bigl(F_{11}(\omega) +
F_{22}(\omega) + 2 F_{12}(\omega )\bigr) \cos(t\omega)
\Biggr\rrvert ^{\alpha} \Biggr]^{1/\alpha},
\]
if we write
\[
\bmm{F}(\omega) = \frac{\partial}{\partial\theta}\Psi(\omega )^*\bm {f}(\omega;
\theta)^{-1} \Big \vert_{\theta=\theta_0} \Psi(\omega).
\]
The confidence intervals for $\theta$ are summarized in the following
Table~\ref{table:5}.
%
\begin{table}[t]
\tabcolsep=0pt
\caption{90\% confidence intervals (and length) for true parameter.
Sample size is $300$ and $\alpha=1.5$}\label{table:5}
\begin{tabular*}{\textwidth}{@{\extracolsep{\fill}}llllll@{}}
\hline
& $b$ & $\theta_0 \approx$ &\multicolumn{2}{l}{E.L.} & (Length) \\
\hline
Case 8 & 0& 0.0000 & $-$0.1685 & 0.1690 & (0.3375) \\
Case 9 &0.3& 0.1755 & \hphantom{$-$}0.0467 & 0.3208 & (0.2741) \\
Case 10 &0.6& 0.3669 & \hphantom{$-$}0.2601 & 0.4920 & (0.2320) \\
Case 11 &0.9& 0.5787 & \hphantom{$-$}0.5046 & 0.6641 & (0.1596) \\ \hline
\end{tabular*}
\end{table}

\begin{table}[b]
\tablewidth=150pt
\tabcolsep=0pt
\caption{Coverage errors of confidence intervals for the parameter
$\theta_0$} \label{tbl:4}
\begin{tabular*}{150pt}{@{\extracolsep{\fill}}lll@{}}
\hline
& \multicolumn{2}{l}{Coverage errors} \\
& \multicolumn{1}{l}{E.L.} & \multicolumn{1}{l}{SAC} \\ \hline
Case 1 & 0.082 & 0.087 \\
Case 2 & 0.089 & 0.096 \\
Case 3 & 0.094 & 0.098 \\
Case 4 & 0.082 & 0.087 \\
Case 5 & 0.053 & 0.056 \\
Case 6 & 0.092 & 0.095 \\
Case 7 & 0.086 & 0.090 \\
Case 8 & 0.011 & --- \\
Case 9 & 0.027 & --- \\
Case 10 & 0.032 & --- \\
Case 11 & 0.049 & --- \\ \hline
\end{tabular*}
\end{table}

We also focus on the one-sided coverage error to evaluate the performances
of the confidence intervals.
Let $\theta^{U}$ and $\theta^{L}$ be the endpoints of a confidence interval.
The one-sided coverage error is given by
\[
\bigl\llvert \Pr\bigl[\bigl\{\theta_0 < \theta^{L}\bigr
\} \cup\bigl\{ \theta^{U} < \theta_0\bigr\} \bigr]-0.1
\bigr\rrvert.
\]
In this time, we calculated the confidence intervals constructed by
both methods for univariate case, and by the empirical likelihood
approach for multivariate case by 1000 times of Monte Carlo simulations.
Namely, we made 1000 confidence intervals $(\theta^L_l,\theta^U_l)$,
$l=1,\ldots,1000$, independently,
and calculate the quantity
\[
\biggl\llvert \frac{\sum^{1000}_{l=1}\mathbb{I} \{\theta_0\notin
 (\theta_l^L,\theta_l^{U} ) \}}{1000} - 0.1 \biggr\rrvert
\]
for each case, where $\mathbb{I}$ denotes the indicator function.
Empirical coverage errors are shown in Table~\ref{tbl:4}.
From this table, the
empirical likelihood confidence intervals are more accurate than those
by the existing method.
Especially, it seems that both methods give
the close coverage probabilities to the nominal level when $\alpha$ is
nearly $2.0$.
On the other hand, we can see that both methods give the close coverage
probabilities to the nominal level
as $n$ increases (case 1, case 6 and case 7).

Furthermore, our results also apply in the multivariate case.
Although the coverage error becomes worse as the pseudo true value gets larger,
it can be seen that the confidence intervals correspondingly becomes
smaller in Table~\ref{table:5}.

\section{Proofs}\label{sec6}
This section provides the proofs of theorems.
The following notation will be used throughout this section.
\[
\bmm{P}_{n}(\bm{\theta}_0) \equiv\frac{1}{n}\sum
^{n}_{t=1}{\bmm{m}}(\lambda_t;
\bm{\theta}_0)\quad\mbox{and}\quad \bmm{S}_{n}(\bm{\theta}_0) \equiv\frac{1}{n}\sum^{n}_{t=1}{
\bmm{m}}(\lambda_t;\bm{\theta}_0){\bmm{m}}(
\lambda_t;\bm{\theta}_0)'.
\]
%
\subsection{Proof for Theorem \texorpdfstring{\protect\ref{theorem:1}}{3.1}}\label{sec6.1}
We start with some auxiliary results.
Recalling (\ref{eq:2.1}), let
\[
\rho_{n,A}(h) = \sum^{n-h}_{t=1}
\tilde{A}_t\tilde{A}_{t+h},\qquad h=1,\ldots,n-1\quad \mbox{and}\quad
T_{n,A}(\omega) = 2\sum^{n-1}_{h=1}
\rho_{n,A}(h)\cos( h\omega).
\]

\begin{lem}\label{lem:1}
\[
\E T_{n,Z}(\omega)=0,\qquad\E T_{n,Z}(\omega)^2\to
\cases{1 & \quad$(\omega\not\equiv0 \bmod\uppi)$,
\cr
2 &\quad $(\omega\equiv0 \bmod\uppi)$,} %
\]
as $n\to\infty$ uniformly in $\alpha\in(0,2]$ and $\sigma>0$.
\end{lem}

\begin{pf}
By symmetry and boundedness of $\tilde{Z}_t$'s, $\E\tilde{Z}_1$ exists
and is equal to $0$.
Furthermore, from the definition of $\tilde{Z}_1,\ldots,\tilde{Z}_n$,
we can see that
$\sum^{n}_{t=1}\tilde{Z}_t^2 = 1$ almost surely, so $\E\tilde{Z}_1^2=1/n$.
Using Chebyshev's inequality, we can see
\[
\Pr \bigl\{ \llvert \tilde{Z}_1\rrvert <\epsilon^{-1/2}n^{-1/2}
\bigr\} > 1-\epsilon
\]
for any $\epsilon>0$. This inequality means $\sqrt{n}\tilde{Z}_1^2$ is
$\mathrm{O}_p(n^{-1/2})$, hence $\sqrt{n}\tilde{Z}_1^2$ converges to $0$ in
distribution uniformly in $\alpha\in(0,2]$. Therefore, by Taylor's
theorem there exists a constant $c$ such that
\begin{eqnarray*}
\E\exp\bigl\{\mathrm{i}\xi\sqrt{n}\tilde{Z}_1^2\bigr\} &=& 1-
\frac{\xi^2}{2} n \E\tilde{Z}_1^4+
\frac{\xi^3\sin(\xi
c)}{6}n^{3/2}\E\tilde{Z}_1^6 + \mathrm{i}
\operatorname{Im} \bigl[\E\exp\bigl\{\mathrm{i}\xi\sqrt {n}\tilde{Z}_1^2
\bigr\} \bigr]
\\
&\to&1
\end{eqnarray*}
uniformly in $\xi\in\R$ by L\'{e}vy's continuity theorem,
where $\operatorname{Im}(z)$ means the imaginary part of a complex number $z$.
So we can conclude $n\E\tilde{Z}_1^4$ converges to $0$ as $n$ tends to
$\infty$.
We also find that $n(n-1)\E\tilde{Z}_1^2\tilde{Z}_2^2$ converges to $1$
by taking expectations on both sides of following identical equation:
%
\begin{eqnarray}
1 = \sum^{n}_{t=1}\tilde{Z}_t^4
+ \sum_{t\neq s}\tilde{Z}_t^2
\tilde {Z}_s^2.\label{eq:6.0}
\end{eqnarray}
Remembering the facts above, let us evaluate the expectations. First,
from symmetry of $\tilde{Z_1}$, it is easy to see that $\E
T_{n,Z}(\omega)$ is exactly
equal to $0$. Next, we expand $T_{n,Z}(\omega)^2$ and obtain that
%
\begin{eqnarray}\label{eq:1}
&&\E T_{n,Z}(\omega)^2
\nonumber\\[-8pt]\\[-8pt]
&&\quad= n(n-1)\E\tilde{Z}_1^2\tilde{Z}_2^2
+ 2 n\E\tilde{Z}_1^2\tilde {Z}_2^2
\sum^{n-1}_{h=1}\cos( 2h\omega) - 2 \E
\tilde{Z}_1^2\tilde {Z}_2^2\sum
^{n-1}_{h=1}h\cos( 2h\omega).\nonumber
\end{eqnarray}
The first term of (\ref{eq:1}) converges to $1$ as $n\to\infty$.
Suppose that $\omega\equiv0 \mod\uppi$, then
\[
2 n\E\tilde{Z}_1^2\tilde{Z}_2^2
\sum^{n-1}_{h=1}\cos(2h\omega) - 2 \E
\tilde{Z}_1^2\tilde{Z}_2^2\sum
^{n-1}_{h=1}h\cos( 2h\omega) = n(n-1)\E
\tilde{Z}_1^2\tilde{Z}_2^2 \to1.
\]
Next, for $\omega\not\equiv0 \bmod\uppi$, the following two identical
equations hold;
\begin{eqnarray*}
\sum^{n-1}_{h=1}\cos( 2h\omega) &=&
\frac{\cos( 2(n-1)\omega) + \cos( 2\omega) -\cos( 2n\omega
)}{2(1-\cos(2\omega))},
\\
\sum^{n-1}_{h=1}h\cos( 2h\omega) &=&
\frac{n\cos( 2(n-1)\omega) - (n-1)\cos( 2n\omega)-1}{2(1-\cos(
2\omega))}.
\end{eqnarray*}
Using these equations, we obtain that
\[
2 n\E\tilde{Z}_1^2\tilde{Z}_2^2
\sum^{n-1}_{h=1}\cos( 2h\omega) - 2 \E
\tilde{Z}_1^2\tilde{Z}_2^2\sum
^{n-1}_{h=1}h\cos( 2h\omega) \to0.
\]
Hence we get desired result.
\end{pf}

\begin{lem}\label{lem:2}
$\sum\sum_{k\neq l} \cov\{ \tilde{I}_{n,Z}(\lambda_k)^2 , \tilde
{I}_{n,Z}(\lambda_l)^2 \} = \mathrm{O}(n)$.
\end{lem}

\begin{pf}
From Brillinger \cite{B:2001},
\[
\cov\bigl\{ \tilde{I}_{n,Z}(\lambda_k)^2 ,
\tilde{I}_{n,Z}(\lambda_l)^2 \bigr\} = \sum
^{8}_{\bm{\nu}:p=1}\prod
^{p}_{j=1}\cum\bigl\{d_{n,Z}(\lambda
_{k_j});k_j\in\nu_j\bigr\},
\]
where the summation is taken over all indecomposable partitions $\bm
{\nu
} = \nu_1\cup\cdots\cup\nu_p$, $p=1,\ldots,8$ of a table
%
\begin{equation}\label{eq:105}
\renewcommand{\arraystretch}{1.2}
\begin{array}{c:c:c:c}
 k\hphantom{k} & \hphantom{k} k \hphantom{k}& \hphantom{k}-k \hphantom{k}&\hphantom{k} -k \\ \hdashline
 l\hphantom{k} & \hphantom{k} l \hphantom{k}& \hphantom{k}-l \hphantom{k}&\hphantom{k} -l
\end{array}
\end{equation}
(see Brillinger \cite{B:2001}), and $d_{n,Z}(\lambda_k) = \sum^{n}_{t=1}\tilde
{Z}_t \exp(\mathrm{i}t\lambda_k)$.
Note that $\cum\{d_{n,Z}(\lambda_{k_1}), \ldots,\allowbreak d_{n,Z}(\lambda
_{k_m})\}$ is 0 for odd $m$.
Let us consider following five partitions;
%
\begin{eqnarray}
\begin{array} {lll} &p=1,\qquad & (k,k,-k,-k,l,l,-l,-l),
\\
&p=2,\qquad & (k,-k,l,-l) \cup(k,-k,l,-l),
\\
&\qquad & (k,-k) \cup(k,-k,l,l,-l,-l),
\\
& \qquad& (l,-l) \cup(k,k,-k,-k,l,-l)
\\
\mbox{and}\quad&p=3,\qquad & (k,-k) \cup(l,-l) \cup(k,-k,l,-l). \end{array}
\label{eq:3}
\end{eqnarray}
First, we show that with different $k$ and $l$ in $\bm{\nu}$,
%
\begin{eqnarray}
\sum_{k\neq l}\sum^{8}_{\bm{\nu}':p=1}
\prod^{p}_{j=1}\cum\bigl\{
d_{n,Z}(\lambda_{k_j});k_j\in\nu_j
\bigr\} = \mathrm{O}(n).\label{eq:103}
\end{eqnarray}
for indecomposable decompositions $\bm{\nu}'=\bm{\nu}\setminus
(\ref{eq:3})$.
However, the proof for (\ref{eq:103}) contains lengthy and complex
algebra, so we confine to giving a representative example here.

Let us consider partitions for $p=4$.
We can evaluate the second order cumulant as
\begin{eqnarray*}
\cum\bigl\{ d_{n,Z}(\lambda_{k}), d_{n,Z}(
\lambda_{l}) \bigr\} &=& \E\tilde{Z}_1^2\sum
^{n}_{t=1}\exp\bigl(\mathrm{i}t(\lambda_k-
\lambda _l)\bigr)
\\
&=& \frac{1}{n}\sum^{n}_{t=1}\exp
\biggl(\mathrm{i}t\frac{2\uppi(k-l)}{n} \biggr)
\\
&=& \cases{ 1 &\quad $(k-l\equiv0\bmod n)$,
\cr
0 &\quad $(k-l\not\equiv0\bmod n)$,}
\end{eqnarray*}
therefore
\begin{eqnarray*}
&&\cum\bigl\{d_{n,Z}(\lambda_{k_1}),d_{n,Z}(
\lambda_{k_2})\bigr\}\cdots\cum\bigl\{ d_{n,Z}(
\lambda_{k_7}),d_{n,Z}(\lambda_{k_8})\bigr\}
\\
&&\quad= \cases{ 1 &\quad $(k_1-k_2, \ldots, k_7-k_8 \equiv0\bmod n)$,
\cr
0 &\quad (otherwise). } %
\end{eqnarray*}
So when $p=4$, we obtain
%
\begin{eqnarray}
\sum_{k\neq l}\prod^{p}_{j=1}
\cum\bigl\{d_{n,Z}(\lambda_{k_j});k_j\in
\nu_j\bigr\} = \mathrm{O}(n)\label{eq:104}
\end{eqnarray}
for any indecomposable partition \eqref{eq:105}.
Similarly, we can check (\ref{eq:104}) for $p=2$ and $3$.
Next, we need to check the cumulants on partitions (\ref{eq:3}).
For simplicity, we introduce generic residual terms $R^{(1)}_n(k,l)$,
$\ldots$, $R^{(4)}_n(k,l)$ such that
$\sum\sum_{k\neq l}R^{(\eta)}_n(k,l)^\gamma= \mathrm{O}(n)$ for $\gamma=1$,$2$,
$\eta=1,2,3$ and 4.
A simple example of $R^{(\eta)}_n(k,l)$ is given as
\[
R^{(\eta)}_n(k,l) = \cases{
\exists(\mbox{constant}) &\quad $(k-l\equiv0\bmod n)$,
\cr
0 &\quad $(k-l\not\equiv0\bmod n)$, } %
\]
and these will appear when we expand the cumulants concerned.
The fourth order joint cumulant on $(\lambda_{k},-\lambda_{k},\lambda
_{l},-\lambda_{l})$ is represented as
%
\begin{eqnarray}\label{eq:6.4}
&&\cum\bigl\{ d_{n,Z}(\lambda_{k}), d_{n,Z}(-
\lambda_{k}), d_{n,Z}(\lambda _{l}),
d_{n,Z}(-\lambda_{l}) \bigr\}
\nonumber\\[-8pt]\\[-8pt]
&&\quad= n\E\tilde{Z}_1^4 + n(n-1) \E\tilde{Z}_1^2
\tilde{Z}_2^2 - 1 + R^{(1)}_n(k,l).
\nonumber
\end{eqnarray}
From (\ref{eq:6.0}),
(\ref{eq:6.4}) becomes
\[
\cum\bigl\{ d_{n,Z}(\lambda_{k}), d_{n,Z}(-
\lambda_{k}), d_{n,Z}(\lambda _{l}),
d_{n,Z}(-\lambda_{l}) \bigr\} = R^{(1)}_n(k,l).
\]

By the same argument as above, and using identical equations
\begin{eqnarray*}
\Bigg\{\sum^{n}_{t=1}\tilde{Z}_t^2
\Biggr\} \Biggl\{\sum^{(\ast)} _{t,s}\tilde
{Z}_t^2\tilde{Z}_s^2 \Biggr\} &=&
\sum^{(\ast)}_{t,s}\tilde{Z}_t^2
\tilde{Z}_s^2 = 2\sum^{(\ast)}_{t,s}
\tilde{Z}_t^4\tilde{Z}_s^2+
\sum^{(\ast
)}_{t,s,u}\tilde {Z}_t^2
\tilde{Z}_s^2\tilde{Z}_u^2,
\\
\Biggl\{\sum^{n}_{t=1}\tilde{Z}_t^2
\Biggr\} \Biggl\{\sum ^{n}_{t=1}
\tilde{Z}_t^4 \Biggr\} &=& \sum^{n}_{t=1}
\tilde{Z}_t^4 = \sum^{(\ast)}_{t,s}
\tilde{Z}_t^4\tilde{Z}_s^2 +
\sum ^{n}_{t=1}\tilde {Z}_t^6
\end{eqnarray*}
and
%
\begin{eqnarray}\label{eq:6.5}
1 &=& \sum^{n}_{t=1}\tilde{Z}_t^8
+ 4\sum^{(\ast)}_{t,s}\tilde {Z}_t^6
\tilde {Z}_s^2 + 3 \sum^{(\ast)}_{t,s}
\tilde{Z}_t^4\tilde {Z}_s^4
\nonumber\\[-8pt]\\[-8pt]
&&{}+ 6\sum^{(\ast)}_{t,s,u}\tilde{Z}_t^4
\tilde{Z}_s^2\tilde{Z}_u^2 +
\sum^{(\ast)} _{t,s,u,v}\tilde{Z}_t^2
\tilde{Z}_s^2\tilde{Z}_u^2\tilde
{Z}_v^2,\nonumber
\end{eqnarray}
we obtain that
%
\begin{eqnarray}\label{eq:4.4}
&&\cum\bigl\{d_{n,Z}(\lambda_{k}),d_{n,Z}(-
\lambda_{k}), d_{n,Z}(\lambda_{l}),
d_{n,Z}(\lambda_{l}),d_{n,Z}(-
\lambda_{l}), d_{n,Z}(-\lambda_{l})
\bigr\}
\nonumber
\\
&&\quad= R^{(2)}_n(k,l),
\nonumber
\\
&&\cum\bigl\{ d_{n,Z}(\lambda_{k}),
\ldots,d_{n,Z}(-\lambda _{l})\bigr\}\qquad\mbox{(the eighth order joint cumulant)}
\\
&&\quad=2 n^2 \E\tilde{Z}_1^4
\tilde{Z}_2^4 - 6 n^3 \E
\tilde{Z}_1^4\tilde{Z}_2^2
\tilde{Z}_3^2
\nonumber
\\
&&\qquad{}+ n^4 \E\tilde{Z}_1^2
\tilde{Z}_2^2\tilde{Z}_3^2
\tilde{Z}_4^2 - \bigl\{ n^2 \E
\tilde{Z}_1^2\tilde{Z}_2^2 \bigr
\}^2 + R^{(3)}_n(k,l),\nonumber
\end{eqnarray}
where ${\sum}^{(\ast)}_{t_1,\ldots,t_m}$ is a summation taken over all
$t_1,\ldots,t_m$ are different from each other.\vspace*{2pt}

According to the same argument as that in Lemma~\ref{lem:1}, the
first and second terms in (\ref{eq:4.4}) converge to $0$ as $n\to
\infty$, and the fourth term converges to $1$.

Finally, from (\ref{eq:6.5}), the third term converges to $1$.
Hence the eighth order joint cumulant becomes
\[
\cum\bigl\{ d_{n,Z}(\lambda_{k}),\ldots,d_{n,Z}(-
\lambda_{l})\bigr\} = R^{(4)}_n(k,l),
\]
so we have\vspace*{-1pt}
\[
\sum_{k\neq l}\sum^{8}_{\bm{\nu}:p=1}
\prod^{p}_{j=1}\cum\bigl\{
d_{n,Z}(\lambda_{k_j});k_j\in\nu_j
\bigr\} = \mathrm{O}(n).
\]
\upqed\end{pf}

\begin{lem}\label{lem:3} Under Assumption~\ref{ass:2.1},\vspace*{-1pt}
\[
\bmm{S}_n(\bm{\theta}_0) \plim\bmm{W}
\]
as $n \to\infty$. Here $\bmm{W}$ is defined in Theorem~3.1.\vspace*{-1pt}
\end{lem}

\begin{pf}
We first make use of the decomposition of the periodogram
in Kl\"uppelberg and Mikosch \cite{KM:1994} as follows, that is,\vspace*{-1pt}
%
\begin{eqnarray}\label{eq:201}
\tilde{I}_{n,X}(\omega)^2 &=& \tilde{g}(
\omega)^2\tilde{I}_{n,Z}(\omega)^2 +
\mathrm{o}_p(1)
\nonumber
\\
&=& \tilde{g}(\omega)^2 \Biggl\{ 1 + 2\sum
^{n-1}_{h=1}\rho _{n,Z}(h)\cos( h\omega)
\Biggr\}^2 + \mathrm{o}_p(1)
\\
&=& \tilde{g}(\omega)^2\bigl\{ 1+2 T_{n,Z}(\omega) +
T_{n,Z}(\omega)^2 \bigr\} + \mathrm{o}_p(1).
\nonumber
\end{eqnarray}
Then from Lemma~\ref{lem:1}, we obtain that\vspace*{-1pt}
\begin{eqnarray*}
\E \bigl[\bmm{S}_n(\bm{\theta}_0) \bigr] &=&
\frac{1}{n}\sum^{n}_{t=1}
\frac{\partial{f(\lambda_t;\bm{\theta})^{-1}}}{\partial\bm
{\theta}} \frac{\partial{f(\lambda_t;\bm{\theta})^{-1}}}{\partial\bm
{\theta}'} \Big|_{\bm{\theta}=\bm{\theta}_0} \E\tilde{I}_{n,X}(
\lambda_t)^2
\\
&\to& \frac{1}{2\uppi}\int^{\uppi}_{-\uppi}
\frac{\partial f(\omega;\bm{\theta
})^{-1}}{\partial\bm{\theta}}\frac{\partial f(\omega;\bm{\theta
})^{-1}}{\partial\bm{\theta}'} \Big|_{\bm{\theta}=\bm{\theta}_0} 2 \tilde{g}(
\omega)^2\,\mathrm{d}\omega= \bmm{W}.
\end{eqnarray*}
From Lemma~\ref{lem:2}, Assumption~\ref{ass:2.1} and \eqref{eq:201}, if
we define\vspace*{-1pt}
\[
h_{\bm{\theta}_0}(\omega)_{ab} = \frac{\partial{f(\omega;\bm{\theta})^{-1}}}{\partial{\theta}_{a}}
\frac{\partial{f(\omega;\bm{\theta})^{-1}}}{\partial{\theta
}_{b}}\Big |_{\bm{\theta}=\bm{\theta}_0}\tilde{g}(\omega)^2,
\]
then\vspace*{-1pt}
\begin{eqnarray*}
\cov\bigl\{ \bmm{S}_n(\bm{\theta}_0)_{ab},
\bmm{S}_n(\bm{\theta }_0)_{cd} \bigr\}
&=& \frac{1}{n^2}\sum^{n}_{t=1}\sum
^{n}_{s=1} h_{\bm{\theta
}_0}(\lambda
_t)_{ab} h_{\bm{\theta}_0}(\lambda_s)_{cd}
\cov\bigl\{ \tilde{I}_{n,Z}(\lambda_t)^2,
\tilde{I}_{n,Z}(\lambda_s)^2 \bigr\}
\\
&=& \frac{1}{n^2}\sum^{n}_{t=1}h_{\bm{\theta}_0}(
\lambda_t)_{ab} h_{\bm
{\theta}_0}(\lambda_t)_{cd}
\var\tilde{I}_{n,Z}(\lambda _t)^2
\\
&&{} +\frac{1}{n^2}\sum_{t\neq s} h_{\bm{\theta}_0}(
\lambda_t)_{ab} h_{\bm{\theta}_0}(\lambda_s)_{cd}
\cov\bigl\{ \tilde{I}_{n,Z}(\lambda_t)^2,
\tilde{I}_{n,Z}(\lambda_s)^2 \bigr\} + \mathrm{o}(1)
\\
&\to& 0
\end{eqnarray*}
for $a$, $b$, $c$, $d=1,\ldots,q$.
These facts imply the convergence of $\bmm{S}_n(\bm{\theta}_0)$ in
probability.
\end{pf}

\begin{pf*}{Proof of Theorem~\ref{theorem:1}}
By Lagrange's multiplier method,
$w_1,\ldots,w_n$ which maximize the objective function in $R(\bm
{\theta
})$ are given by
\[
w_t = \frac{1}{n}\frac{1}{1+\phi'{\bmm{m}}(\lambda_t;\bm{\theta
}_0)},\qquad t=1,\ldots, n,
\]
where $\phi\in\R^q$ is
the Lagrange multiplier which is defined as the solution of $q$-restrictions
%
\begin{equation}\label{eq:101}
J_{n,\bm{\theta}_0}(\phi) = \frac{1}{n}\sum^{n}_{t=1}
\frac{{\bmm{m}}(\lambda_t;\bm{\theta}_0)}{1+\phi'{\bmm{m}}(\lambda_t;\bm
{\theta
}_0)} = \bm{0}.
\end{equation}
First of all, let us derive the order of $\phi$.
Set $Y_t \equiv\phi' {\bmm{m}}(\lambda_t;\bm{\theta}_0)$ and from
(\ref{eq:101}),
\begin{eqnarray*}
\bm{0} &=& \frac{1}{n}\sum^{n}_{t=1}
\frac{{\bmm{m}}(\lambda_t;\bm{\theta
}_0)}{1+Y_t}
\\
&=& \frac{1}{n}\sum^{n}_{t=1}
\biggl\{1-Y_t+\frac{Y_t^2}{1+Y_t} \biggr\} {\bmm{m}}(
\lambda_t;\bm{\theta}_0)
\\
&=& \bmm{P}_n(\bm{\theta}_0)-\bmm{S}_n(\bm{
\theta}_0)\phi+\frac
{1}{n}\sum^{n}_{t=1}
\frac{{\bmm{m}}(\lambda_t;\bm{\theta
}_0)Y_t^2}{1+Y_t}.
\end{eqnarray*}
Hence,
%
\begin{equation}
\label{eq:102} \phi= \bmm{S}_n(\bm{\theta}_0)^{-1}
\Biggl\{\bmm{P}_n(\bm{\theta }_0)+\frac
{1}{n}\sum
^{n}_{t=1}\frac{{\bmm{m}}(\lambda_t;\bm{\theta
}_0)Y_t^2}{1+Y_t} \Biggr\}
\equiv \bmm{S}_n(\bm{\theta}_0)^{-1}
\bmm{P}_n(\bm{\theta}_0) + \bm{\epsilon}\qquad\mbox{(say)}.
\end{equation}
Next, we introduce $M_n\equiv\max_{1\leq k \leq n}\llVert {\bmm{m}}(\lambda
_k;\bm{\theta}_0)\rrVert _E$.
The order of $M_n$ is given by
\begin{eqnarray*}
M_n &=& \max_{1\leq t \leq n} \biggl\llVert
\frac{\partial f(\lambda_t;\bm
{\theta
})^{-1}}{\partial\bm{\theta}}\Big |_{\bm{\theta}=\bm{\theta
}_0}\tilde {I}_{n,X}(\lambda_t)
\biggr\rrVert _E
\\
&\leq&\max_{1\leq t \leq n} \biggl\llVert \frac{\partial f(\lambda_t;\bm
{\theta
})^{-1}}{\partial\bm{\theta}}
\Big|_{\bm{\theta}=\bm{\theta
}_0} \biggr\rrVert _E\max_{1\leq t \leq n}\bigl
\llvert I_{n,X}(\lambda_t)\bigr\rrvert \frac{1}{\gamma
_{n,X}^2}
\\
&\leq&\max_{\omega\in[-\uppi,\uppi]} \biggl\llVert \frac{\partial f(\omega;{\bm
{\theta}})^{-1}}{\partial\bm{\theta}}
\Big|_{\bm{\theta}=\bm
{\theta
}_0} \biggr\rrVert _E\max_{\omega\in[-\uppi,\uppi]}\bigl
\llvert I_{n,X}(\omega)\bigr\rrvert \frac
{1}{\gamma
_{n,X}^2}
\\
&=& \max_{\omega\in[-\uppi,\uppi]} \biggl\llVert \frac{\partial f(\omega;{\bm{\theta
}})^{-1}}{\partial\bm{\theta}}
\Big|_{\bm{\theta}=\bm{\theta
}_0} \biggr\rrVert _E \max_{\omega\in[-\uppi,\uppi]}\bigl
\llvert g(\omega)\bigr\rrvert \frac{\max_{\omega\in
[-\uppi,\uppi
]}\llvert I_{n,X}(\omega)\rrvert }{\max_{\omega\in[-\uppi,\uppi]}\llvert g(\omega)\rrvert }\frac
{1}{\gamma_{n,X}^2}
\\
&\leq&\max_{\omega\in[-\uppi,\uppi]} \biggl\llVert \frac{\partial f(\omega;{\bm
{\theta}})^{-1}}{\partial\bm{\theta}}
\Big|_{\bm{\theta}=\bm
{\theta
}_0} \biggr\rrVert _E \max_{\omega\in[-\uppi,\uppi]}\bigl
\llvert g(\omega)\bigr\rrvert \max_{\omega\in[-\uppi,\uppi
]} \biggl\llvert
\frac{I_{n,X}(\omega)}{g(\omega)} \biggr\rrvert \frac{1}{\gamma
_{n,X}^2}
\\
&=&{\exists}c_0\max_{\omega\in[-\uppi,\uppi]} \biggl\llvert
\frac
{I_{n,X}(\omega
)}{g(\omega)} \biggr\rrvert\qquad(\because\mbox{ Assumption }\ref{ass:2.1}).
\end{eqnarray*}
On the other hand, it is not difficult to check that
Assumption~\ref{ass:2.2} is sufficient condition for
Corollary~3.3 of Mikosch, Resnick and Samorodnitsky \cite{KSG:2000},
so we have
$M_n=\mathrm{O}_p(\beta_n^2)$, where
\[
\beta_n=
\cases{ (\log n)^{1-1/\alpha} &\quad
$(1<\alpha<2)$,
\cr
\log\log n &\quad $(\alpha=1)$.} %
\]

Henceforth, let $1<\alpha<2$. In the case of $\alpha=1$, the same
argument as follows will go on.
By Ogata and Taniguchi \cite{OT:2010}, there exists a unit vector $\bmm{u}$ in
$\R^q$ such that the following inequality holds:
\[
\Vert\phi\Vert_E \bigl\{ \bmm{u}'
\bmm{S}_n(\bm{\theta}_0)\bmm{u}-\bmm{u}'M_n
\bmm{P}_n(\bm{\theta}_0) \bigr\} \leq
\bmm{u}'\bmm{P}_n(\bm {\theta }_0).
\]
Lemma P5.1 of Brillinger \cite{B:2001} allows us to write $x_n \bmm{P_n}(\bm
{\theta}_0)$ as
\begin{eqnarray*}
x_n\bmm{P}_n(\bm{\theta}_0) &=&
\frac{1}{2\uppi}x_n\int^{\uppi}_{-\uppi}
\frac{\partial f(\omega;\bm
{\theta
})}{\partial\bm{\theta}} \Big|_{\bm{\theta}=\bm{\theta
}_0}\tilde {I}_{n,X}(\omega)\,\mathrm{d}\omega+
\mathrm{O}_p \biggl( \frac{x_n}{n} \biggr)
\\
&=& \frac{1}{2\uppi}\frac{1}{\gamma^2_{n,X}} x_n \int
^{\uppi}_{-\uppi
}\frac
{\partial
f(\omega;\bm{\theta})}{\partial\bm{\theta}} \Big|_{\bm{\theta
}=\bm
{\theta}_0}
\bigl\{{I}_{n,X}(\omega) - T_n \psi^2 \tilde{g}(
\omega) \bigr\}\,\mathrm{d}\omega + \mathrm{O}_p \biggl(\frac{x_n}{n} \biggr),
\end{eqnarray*}
where
\[
T_n = \frac{1}{2\uppi}\int^{\uppi}_{-\uppi}
\frac{I_{n,X}(\omega)}{\psi
^2 \tilde
{g}(\omega)}\,\mathrm{d}\omega.
\]
Then, by Proposition~3.5 of Kl\"{u}ppelberg and Mikosch \cite{KM:1996}
and Cram\'er--Wold device, we have
\begin{eqnarray*}
&&
\pmatrix{ \gamma_{n,X}^2 \cr
\displaystyle
x_n\int^{\uppi}_{-\uppi}\frac{\partial f(\omega;\bm{\theta
})}{\partial{\theta
_1}}
\Big|_{\bm{\theta}=\bm{\theta}_0} \bigl\{{I}_{n,X}(\omega) - T_n \psi
^2 \tilde{g}(\omega) \bigr\}\,\mathrm{d}\omega \cr
\vdots \cr
\displaystyle
x_n\int^{\uppi}_{-\uppi}\frac{\partial f(\omega;\bm{\theta
})}{\partial{\theta
_q}}
\Big|_{\bm{\theta}=\bm{\theta}_0} \bigl\{{I}_{n,X}(\omega) - T_n \psi
^2 \tilde{g}(\omega) \bigr\}\,\mathrm{d}\omega }
\cr
&&\quad\dlim
\pmatrix{
 \psi^2
S_0
\cr
\displaystyle 2\sum^{\infty}_{t=1} S_t \biggl\{
\int^{\uppi}_{-\uppi}\frac{\partial
f(\omega;\bm
{\theta})}{\partial{\theta_1}} \Big|_{\bm{\theta}=\bm{\theta
}_0}
\psi^2 \tilde{g}(\omega)\cos(t\omega)\,\mathrm{d}\omega \biggr\}
\cr
\vdots
\cr
\displaystyle 2\sum^{\infty}_{t=1} S_t \biggl\{
\int^{\uppi}_{-\uppi}\frac{\partial
f(\omega;\bm
{\theta})}{\partial{\theta_q}}\Big |_{\bm{\theta}=\bm{\theta
}_0}
\psi^2 \tilde{g}(\omega)\cos(t\omega)\,\mathrm{d}\omega \biggr\} }.
\end{eqnarray*}
Therefore
%
\begin{eqnarray}
x_n \bmm{P}_n(\bm{\theta}_0)\dlim
\bmm{V}\label{eq:100}
\end{eqnarray}
for $\alpha\in[1,2)$ as $n\to\infty$, where $\bmm{V}$ is defined in
Theorem~\ref{theorem:1}.
So we obtain
%
\begin{equation}\label{eq:6.9}
\mathrm{O}_p\bigl(\Vert\phi\Vert_E\bigr) \bigl[ \mathrm{O}_p(1)-
\underline{\mathrm{O}_p \bigl\{ (\log n)^{2-2/\alpha} \bigr\}\cdot
\mathrm{O}_p\bigl(x_n^{-1}\bigr)} \bigr] \leq
\mathrm{O}_p\bigl(x_n^{-1}\bigr).
\end{equation}
Because as $n\to\infty$,
\begin{eqnarray*}
(\log n)^{2-2/\alpha}x_n^{-1} &=& (\log
n)^{2-2/\alpha} \biggl(\frac{\log n}{n} \biggr)^{1/\alpha
}
\\
&=& \frac{1}{(\log n)^{1/\alpha}} \frac{(\log n)^2}{n^{1/\alpha
}}
\\
&\to& 0,
\end{eqnarray*}
the underlined part in (\ref{eq:6.9}) is $\mathrm{O}_p(1)$. Therefore, we obtain
%
\begin{equation}\label{eq:110}
\mathrm{O}_p\bigl(\Vert\phi\Vert_E\bigr) \leq \mathrm{O}_p
\bigl(x_n^{-1}\bigr).
\end{equation}
On the other hand,
%
\begin{eqnarray}\label{eq:111}
\frac{1}{n}\sum^{n}_{t=1}\bigl\Vert{
\bmm{m}}(\lambda_t;\bm{\theta }_0)\bigr\Vert_E^3
&=& \frac{1}{n}\sum^{n}_{t=1}\bigl\Vert{
\bmm{m}}(\lambda_t;\bm{\theta }_0)\bigr\Vert_E
\bigl\Vert{\bmm{m}}(\lambda_t;\bm{\theta_0})\bigr\Vert
_E^2
\nonumber
\\
&\leq& \frac{1}{n}\sum^{n}_{t=1}M_n
{\bmm{m}}(\lambda_t;\bm{\theta}_0)' {
\bmm{m}}(\lambda_t;\bm{\theta}_0)
\nonumber\\[-8pt]\\[-8pt]
&=& M_n \tr \bigl\{\bmm{S}_n(\bm{\theta}_0)
\bigr\}
\nonumber
\\
&=& \mathrm{O}_p \bigl\{ (\log n)^{2-2/\alpha} \bigr\}.\nonumber
\end{eqnarray}
From (\ref{eq:110}) and (\ref{eq:111}), $\bm{\epsilon}$ in
(\ref{eq:102}) satisfies
%
\begin{equation}
\Vert\bm{\epsilon}\Vert_E \leq\frac{1}{n}\sum
^{n}_{t=1}\bigl\Vert{\bmm{m}}(\lambda_t;\bm{
\theta })\bigr\Vert _E^3 \Vert\phi\Vert_E^2
\llvert 1+Y_t\rrvert ^{-1}.
\end{equation}
Thus, we have
\[
 \mathrm{O}_p\bigl(\Vert x_n \bm{\epsilon}\Vert_E\bigr) =
\mathrm{O}_p \biggl\{ \frac{(\log
n)^{2-1/\alpha}}{n^{1/\alpha}} \biggr\} \plim 0.
\]

Now let us show the convergence of the empirical likelihood ratio statistic.
Under $H$: $\bm{\theta}=\bm{\theta}_0$, $-2({x_n^2}/{n})\log
{R}(\bm
{\theta}_0)$ can be expanded as
\begin{eqnarray*}
-2\frac{x_n^2}{n}\log{R}(\bm{\theta}_0) &=& -2\frac{x_n^2}{n}
\sum^{n}_{t=1}\log n w_t
\\
&=& 2\frac{x_n^2}{n}\sum^{n}_{t=1}
\log(1+Y_t)
\\
&=& 2\frac{x_n^2}{n}\sum^{n}_{t=1}Y_t-
\frac{x_n^2}{n}\sum^{n}_{t=1}Y_t^2+2
\frac{x_n^2}{n}\sum^{n}_{t=1}\mathrm{O}_p
\bigl(Y_t^3\bigr),
\end{eqnarray*}
where
\begin{eqnarray*}
2\frac{x_n^2}{n}\sum^{n}_{t=1}Y_t
&=& 2\frac{x_n^2}{n}\sum^{n}_{t=1}
\phi' {\bmm{m}}(\lambda_t;\bm {\theta }_0)
\\
&=& 2\frac{x_n^2}{n} \bigl\{ {\bmm{S}}_n(\bm{
\theta}_0)^{-1}{\bmm{P}}_n(\bm{
\theta}_0)+{\bm{\epsilon}} \bigr\}' \sum
^{n}_{t=1}{\bmm{m}}(\lambda_t;\bm{
\theta}_0)
\\
&=& 2 x_n^2 \bigl\{ {\bmm{P}}_n(\bm{
\theta}_0)'{\bmm{S}}_n(\bm {\theta
}_0)^{-1}+{\bm{\epsilon}}' \bigr\} {
\bmm{P}}_n(\bm{\theta}_0)
\\
&=& 2 \bigl\{ x_n {\bmm{P}}_n(\bm{\theta}_0)
\bigr\}' {\bm {S}}_n(\bm {\theta}_0)^{-1}
\bigl\{ x_n {\bmm{P}}_n(\bm{\theta}_0) \bigr
\} + 2 (x_n{\bm{\epsilon}})' \bigl\{x_n{
\bmm{P}}_n(\bm{\theta }_0) \bigr\},
\\
\frac{x_n^2}{n}\sum^{n}_{t=1}Y_t^2
&=& \frac{x_n^2}{n}\sum^{n}_{t=1} \bigl
\{ \phi'{\bmm{m}}(\lambda _t;\bm {\theta}_0)
\bigr\}^2
\\
&=& x_n^2 \phi' {\bmm{S}}_n(
\bm{\theta}_0) \phi
\\
&=& x_n^2 \bigl\{ {\bmm{P}}_n(\bm{
\theta}_0)'{\bmm{S}}_n(\bm{\theta
}_0)^{-1} + {\bm{\epsilon}}' \bigr\} {
\bmm{S}}_n(\bm{\theta}_0) \bigl\{ {\bmm{S}}_n(
\bm{\theta}_0)^{-1}{\bmm{P}}_n(\bm{
\theta}_0) + {\bm {\epsilon}} \bigr\}
\\
&=& \bigl\{ x_n {\bmm{P}}_n(\bm{\theta}_0)
\bigr\}' {\bmm{S}}_n(\bm {\theta}_0)^{-1}
\bigl\{ x_n {\bmm{P}}_n(\bm{\theta}_0) \bigr
\}
\\
& &{} + ( x_n {\bm{\epsilon}} )' {\bmm{S}}_n(
\bm{\theta }_0) ( x_n {\bm{\epsilon}} ) + 2 (
x_n {\bm{\epsilon}} )' \bigl\{ x_n {
\bmm{P}}_n(\bm {\theta }_0) \bigr\}
\end{eqnarray*}
and
\begin{eqnarray*}
\frac{x_n^2}{n} \Biggl\llvert \sum^{n}_{t=1}\mathrm{O}_p
\bigl(Y_t^3\bigr) \Biggr\rrvert &\leq&
\frac{x_n^2}{n}\exists c\sum^{n}_{t=1}
\llvert Y_t\rrvert ^3
\\
&=& \frac{x_n^2}{n}c\Vert\phi\Vert_E^3 \sum
^{n}_{t=1}\bigl\Vert{\bmm{m}}(
\lambda_t;\bm{\theta}_0)\bigr\Vert_E^3
\\
&=& \frac{x_n^2}{n}\mathrm{O}_p\bigl(x_n^{-3}\bigr)
\cdot \mathrm{O}_p\bigl\{ n(\log n)^{2-2/\alpha}\bigr\}
\\
&=& \mathrm{O}_p \biggl\{ \frac{(\log n)^{2-1/\alpha}}{n^{1/\alpha}} \biggr\}
\\
&\plim & 0\qquad(n\to\infty).
\end{eqnarray*}
Hence, using (\ref{eq:100}) and Lemma~\ref{lem:3},
\begin{eqnarray*}
-\frac{2x_n^2}{n}\log{R}(\bm{\theta}_0)& =& \bigl\{
x_n {\bmm{P}}_n(\bm{\theta}_0) \bigr
\}' {\bmm{S}}_n(\bm {\theta}_0)^{-1}
\bigl\{ x_n {\bmm{P}}_n(\bm{\theta}_0) \bigr
\} + \mathrm{o}_p(1)
\\
&\dlim & \bmm{V}' \bmm{W}^{-1} \bmm{V}
\end{eqnarray*}
for $\alpha\in[1,2)$.
\end{pf*}


\section*{Acknowledgments}
We thank the area editor, the associate editor and the referees for
their constructive comments
on an earlier draft of the paper. M. Taniguchi was supported by the
Japanese Grant-in-Aid:
A1150300.

\begin{supplement}
\stitle{Proof of Theorem~\ref{theorem:4.1}}
\slink[doi]{10.3150/14-BEJ636SUPP} 
\sdatatype{.pdf}
\sfilename{BEJ636\_supp.pdf}
\sdescription{We provide additional supporting material for the proof of Theorem~\ref{theorem:4.1}.}
\end{supplement}


\printhistory

\begin{thebibliography}{30}


\bibitem{akashi}
\begin{barticle}[author]
\bauthor{\bsnm{Akashi},~\bfnm{Fumiya}\binits{F.}},
\bauthor{\bsnm{Liu},~\bfnm{Yan}\binits{Y.}} \AND
\bauthor{\bsnm{Taniguchi},~\bfnm{Masanobu}\binits{M.}}
(\byear{2014}).
\btitle{Supplement to ``An empirical likelihood approach for symmetric $\alpha $-stable processes.''
DOI:\doiurl{10.3150/14-BEJ636SUPP}}.
\end{barticle}
\bptok{imsref}%
\endbibitem


\bibitem{bhansali1980autoregressive}
\begin{barticle}[mr]
\bauthor{\bsnm{Bhansali},~\bfnm{R.~J.}\binits{R.J.}}
(\byear{1980}).
\btitle{Autoregressive and window estimates of the inverse correlation function}.
\bjournal{Biometrika}
\bvolume{67}
\bpages{551--566}.
\bid{doi={10.1093/biomet/67.3.551}, issn={0006-3444}, mr={0601091}}
\end{barticle}
\bptok{imsref}%
\endbibitem

\bibitem{B:2001}
\begin{bbook}[mr]
\bauthor{\bsnm{Brillinger},~\bfnm{David~R.}\binits{D.R.}}
(\byear{2001}).
\btitle{Time Series}.
\bseries{Classics in Applied Mathematics}
\bvolume{36}.
\blocation{Philadelphia, PA}:
\bpublisher{Society for Industrial and Applied Mathematics (SIAM)}.
\bnote{Data analysis and theory, Reprint of the 1981 edition}.
\bid{doi={10.1137/1.9780898719246}, mr={1853554}}
\end{bbook}
\bptok{imsref}%
\endbibitem

\bibitem{BD:1991}
\begin{bbook}[mr]
\bauthor{\bsnm{Brockwell},~\bfnm{Peter~J.}\binits{P.J.}} \AND
\bauthor{\bsnm{Davis},~\bfnm{Richard~A.}\binits{R.A.}}
(\byear{1991}).
\btitle{Time Series: Theory and Methods},
\bedition{2nd} ed.
\bseries{Springer Series in Statistics}.
\blocation{New York}:
\bpublisher{Springer}.
\bid{doi={10.1007/978-1-4419-0320-4}, mr={1093459}}
\end{bbook}
\bptok{imsref}%
\endbibitem

\bibitem{cleveland1972inverse}
\begin{barticle}[author]
\bauthor{\bsnm{Cleveland},~\bfnm{William~S.}\binits{W.S.}}
(\byear{1972}).
\btitle{The inverse autocorrelations of a time series and their applications}.
\bjournal{Technometrics}
\bvolume{14}
\bpages{277--293}.
\end{barticle}
\bptok{imsref}%
\endbibitem

\bibitem{DR1985a}
\begin{barticle}[mr]
\bauthor{\bsnm{Davis},~\bfnm{Richard}\binits{R.}} \AND
\bauthor{\bsnm{Resnick},~\bfnm{Sidney}\binits{S.}}
(\byear{1985}).
\btitle{Limit theory for moving averages of random variables with regularly varying tail probabilities}.
\bjournal{Ann. Probab.}
\bvolume{13}
\bpages{179--195}.
\bid{issn={0091-1798}, mr={0770636}}
\end{barticle}
\bptok{imsref}%
\endbibitem

\bibitem{DR1985b}
\begin{barticle}[mr]
\bauthor{\bsnm{Davis},~\bfnm{Richard}\binits{R.}} \AND
\bauthor{\bsnm{Resnick},~\bfnm{Sidney}\binits{S.}}
(\byear{1985}).
\btitle{More limit theory for the sample correlation function of moving averages}.
\bjournal{Stochastic Process. Appl.}
\bvolume{20}
\bpages{257--279}.
\bid{doi={10.1016/0304-4149(85)90214-5}, issn={0304-4149}, mr={0808161}}
\end{barticle}
\bptok{imsref}%
\endbibitem

\bibitem{DR1986}
\begin{barticle}[mr]
\bauthor{\bsnm{Davis},~\bfnm{Richard}\binits{R.}} \AND
\bauthor{\bsnm{Resnick},~\bfnm{Sidney}\binits{S.}}
(\byear{1986}).
\btitle{Limit theory for the sample covariance and correlation functions of moving averages}.
\bjournal{Ann. Statist.}
\bvolume{14}
\bpages{533--558}.
\bid{doi={10.1214/aos/1176349937}, issn={0090-5364}, mr={0840513}}
\end{barticle}
\bptok{imsref}%
\endbibitem


\bibitem{DHR:2000}
\begin{barticle}[mr]
\bauthor{\bsnm{Drees},~\bfnm{Holger}\binits{H.}},
\bauthor{\bparticle{de} \bsnm{Haan},~\bfnm{Laurens}\binits{L.}} \AND
\bauthor{\bsnm{Resnick},~\bfnm{Sidney}\binits{S.}}
(\byear{2000}).
\btitle{How to make a {H}ill plot}.
\bjournal{Ann. Statist.}
\bvolume{28}
\bpages{254--274}.
\bid{doi={10.1214/aos/1016120372}, issn={0090-5364}, mr={1762911}}
\end{barticle}
\bptok{imsref}%
\endbibitem

\bibitem{Fama1965}
\begin{barticle}[author]
\bauthor{\bsnm{Fama},~\bfnm{Eugene~F.}\binits{E.F.}}
(\byear{1965}).
\btitle{The behavior of stock-market prices}.
\bjournal{J. Bus.}
\bvolume{38}
\bpages{34--105}.
\end{barticle}
\bptok{imsref}%
\endbibitem

\bibitem{label2568}
\begin{barticle}[mr]
\bauthor{\bsnm{Hall},~\bfnm{Peter}\binits{P.}}
(\byear{1982}).
\btitle{On some simple estimates of an exponent of regular variation}.
\bjournal{J. R. Stat. Soc. Ser. B Stat. Methodol.}
\bvolume{44}
\bpages{37--42}.
\bid{issn={0035-9246}, mr={0655370}}
\end{barticle}
\bptok{imsref}%
\endbibitem

\bibitem{hannan2009multiple}
\begin{bbook}[mr]
\bauthor{\bsnm{Hannan},~\bfnm{E.~J.}\binits{E.J.}}
(\byear{1970}).
\btitle{Multiple Time Series}.
\blocation{New York}:
\bpublisher{Wiley}.
\bid{mr={0279952}}
\end{bbook}
\bptok{imsref}%
\endbibitem

\bibitem{H:1991}
\begin{barticle}[mr]
\bauthor{\bsnm{Hsing},~\bfnm{Tailen}\binits{T.}}
(\byear{1991}).
\btitle{On tail index estimation using dependent data}.
\bjournal{Ann. Statist.}
\bvolume{19}
\bpages{1547--1569}.
\bid{doi={10.1214/aos/1176348261}, issn={0090-5364}, mr={1126337}}
\end{barticle}
\bptok{imsref}%
\endbibitem

\bibitem{KM:1993}
\begin{barticle}[mr]
\bauthor{\bsnm{Kl{\"u}ppelberg},~\bfnm{Claudia}\binits{C.}} \AND
\bauthor{\bsnm{Mikosch},~\bfnm{Thomas}\binits{T.}}
(\byear{1993}).
\btitle{Spectral estimates and stable processes}.
\bjournal{Stochastic Process. Appl.}
\bvolume{47}
\bpages{323--344}.
\bid{doi={10.1016/0304-4149(93)90021-U}, issn={0304-4149}, mr={1239844}}
\end{barticle}
\bptok{imsref}%
\endbibitem

\bibitem{KM:1994}
\begin{barticle}[mr]
\bauthor{\bsnm{Kl{\"u}ppelberg},~\bfnm{Claudia}\binits{C.}} \AND
\bauthor{\bsnm{Mikosch},~\bfnm{Thomas}\binits{T.}}
(\byear{1994}).
\btitle{Some limit theory for the self-normalised periodogram of stable processes}.
\bjournal{Scand. J. Stat.}
\bvolume{21}
\bpages{485--491}.
\bid{issn={0303-6898}, mr={1310091}}
\end{barticle}
\bptok{imsref}%
\endbibitem

\bibitem{KM:1996}
\begin{barticle}[mr]
\bauthor{\bsnm{Kl{\"u}ppelberg},~\bfnm{Claudia}\binits{C.}} \AND
\bauthor{\bsnm{Mikosch},~\bfnm{Thomas}\binits{T.}}
(\byear{1996}).
\btitle{The integrated periodogram for stable processes}.
\bjournal{Ann. Statist.}
\bvolume{24}
\bpages{1855--1879}.
\bid{doi={10.1214/aos/1069362301}, issn={0090-5364}, mr={1421152}}
\end{barticle}
\bptok{imsref}%
\endbibitem

\bibitem{M1963}
\begin{barticle}[author]
\bauthor{\bsnm{Mandelbrot},~\bfnm{Benoit~B.}\binits{B.B.}}
(\byear{1963}).
\btitle{New methods in statistical economics}.
\bjournal{ J. Polit. Econ.}
\bvolume{71}
\bpages{421--440}.
\end{barticle}
\bptok{imsref}%
\endbibitem

\bibitem{MGKA:1995}
\begin{barticle}[mr]
\bauthor{\bsnm{Mikosch},~\bfnm{Thomas}\binits{T.}},
\bauthor{\bsnm{Gadrich},~\bfnm{Tamar}\binits{T.}},
\bauthor{\bsnm{Kl{\"u}ppelberg},~\bfnm{Claudia}\binits{C.}} \AND
\bauthor{\bsnm{Adler},~\bfnm{Robert~J.}\binits{R.J.}}
(\byear{1995}).
\btitle{Parameter estimation for {ARMA} models with infinite variance innovations}.
\bjournal{Ann. Statist.}
\bvolume{23}
\bpages{305--326}.
\bid{doi={10.1214/aos/1176324469}, issn={0090-5364}, mr={1331670}}
\end{barticle}
\bptok{imsref}%
\endbibitem

\bibitem{KSG:2000}
\begin{barticle}[mr]
\bauthor{\bsnm{Mikosch},~\bfnm{Thomas}\binits{T.}},
\bauthor{\bsnm{Resnick},~\bfnm{Sidney}\binits{S.}} \AND
\bauthor{\bsnm{Samorodnitsky},~\bfnm{Gennady}\binits{G.}}
(\byear{2000}).
\btitle{The maximum of the periodogram for a heavy-tailed sequence}.
\bjournal{Ann. Probab.}
\bvolume{28}
\bpages{885--908}.
\bid{doi={10.1214/aop/1019160264}, issn={0091-1798}, mr={1782277}}
\end{barticle}
\bptok{imsref}%
\endbibitem

\bibitem{M:1997}
\begin{barticle}[mr]
\bauthor{\bsnm{Monti},~\bfnm{Anna~Clara}\binits{A.C.}}
(\byear{1997}).
\btitle{Empirical likelihood confidence regions in time series models}.
\bjournal{Biometrika}
\bvolume{84}
\bpages{395--405}.
\bid{doi={10.1093/biomet/84.2.395}, issn={0006-3444}, mr={1467055}}
\end{barticle}
\bptok{imsref}%
\endbibitem

\bibitem{N:2012}
\begin{bbook}[author]
\bauthor{\bsnm{Nolan},~\bfnm{John~P.}\binits{J.P.}}
(\byear{2015}).
\btitle{Stable Distributions -- Models for Heavy Tailed Data}.
\blocation{Boston}:
\bpublisher{Birkh\"{a}user}. To appear.
\end{bbook}
\bptok{imsref}%
\endbibitem

\bibitem{OT:2010}
\begin{barticle}[mr]
\bauthor{\bsnm{Ogata},~\bfnm{Hiroaki}\binits{H.}} \AND
\bauthor{\bsnm{Taniguchi},~\bfnm{Masanobu}\binits{M.}}
(\byear{2010}).
\btitle{An empirical likelihood approach for non-{G}aussian vector stationary processes and its application to minimum contrast estimation}.
\bjournal{Aust. N. Z. J. Stat.}
\bvolume{52}
\bpages{451--468}.
\bid{doi={10.1111/j.1467-842X.2010.00585.x}, issn={1369-1473}, mr={2791530}}
\end{barticle}
\bptok{imsref}%
\endbibitem

\bibitem{O:1988}
\begin{barticle}[mr]
\bauthor{\bsnm{Owen},~\bfnm{Art~B.}\binits{A.B.}}
(\byear{1988}).
\btitle{Empirical likelihood ratio confidence intervals for a single functional}.
\bjournal{Biometrika}
\bvolume{75}
\bpages{237--249}.
\bid{doi={10.1093/biomet/75.2.237}, issn={0006-3444}, mr={0946049}}
\end{barticle}
\bptok{imsref}%
\endbibitem

\bibitem{P:1975}
\begin{bbook}[mr]
\bauthor{\bsnm{Petrov},~\bfnm{V.~V.}\binits{V.V.}}
(\byear{1975}).
\btitle{Sums of Independent Random Variables}.
\blocation{New York}:
\bpublisher{Springer}.
\bid{mr={0388499}}
\end{bbook}
\bptok{imsref}%
\endbibitem

\bibitem{RS:1998}
\begin{barticle}[mr]
\bauthor{\bsnm{Resnick},~\bfnm{Sidney}\binits{S.}} \AND
\bauthor{\bsnm{St{\u{a}}ric{\u{a}}},~\bfnm{Catalin}\binits{C.}}
(\byear{1998}).
\btitle{Tail index estimation for dependent data}.
\bjournal{Ann. Appl. Probab.}
\bvolume{8}
\bpages{1156--1183}.
\bid{doi={10.1214/aoap/1028903376}, issn={1050-5164}, mr={1661160}}
\end{barticle}
\bptok{imsref}%
\endbibitem

\bibitem{RS:1996}
\begin{barticle}[author]
\bauthor{\bsnm{Resnick},~\bfnm{Sidney~I.}\binits{S.I.}} \AND
\bauthor{\bsnm{St{\u{a}}ric{\u{a}}},~\bfnm{Catalin}\binits{C.}}
(\byear{1996}).
\btitle{Asymptotic behavior of Hill's estimator for autoregressive data}.
\bjournal{Stoch. Models}
\bvolume{13}
\bpages{703--723}.
\end{barticle}
\bptok{imsref}%
\endbibitem

\bibitem{rosinski1987multilinear}
\begin{barticle}[mr]
\bauthor{\bsnm{Rosi{\'n}ski},~\bfnm{Jan}\binits{J.}} \AND
\bauthor{\bsnm{Woyczy{\'n}ski},~\bfnm{Wojbor~A.}\binits{W.A.}}
(\byear{1987}).
\btitle{Multilinear forms in {P}areto-like random variables and product random measures}.
\bjournal{Colloq. Math.}
\bvolume{51}
\bpages{303--313}.
\bid{issn={0010-1354}, mr={0891300}}
\end{barticle}
\bptok{imsref}%
\endbibitem

\bibitem{samoradnitsky1994stable}
\begin{bbook}[mr]
\bauthor{\bsnm{Samorodnitsky},~\bfnm{Gennady}\binits{G.}} \AND
\bauthor{\bsnm{Taqqu},~\bfnm{Murad~S.}\binits{M.S.}}
(\byear{1994}).
\btitle{Stable Non-{G}aussian Random Processes}.
\bseries{Stochastic Modeling. Stochastic Models with Infinite Variance}.
\blocation{New York}:
\bpublisher{Chapman and Hall}.
\bid{mr={1280932}}
\end{bbook}
\bptok{imsref}%
\endbibitem

\end{thebibliography}
\end{document}